\documentclass[12pt]{amsart}
\usepackage{amssymb}
\usepackage{amsmath}
\usepackage{longtable}
\newcommand\g{{\mathfrak g}}
\newcommand\h{{\mathfrak h}}
\newcommand{\q}{\mathfrak{q}}
\newcommand{\lf}{\mathfrak{l}}

\newcommand{\Dcal}{\mathcal{D}}
\newcommand{\Fl}{\mathcal{B}}

\newcommand\m{\mathfrak m}
\newcommand\n{\mathfrak n}
\newcommand\borel{\mathfrak b}
\newcommand\z{\mathfrak z}

\renewcommand\t{\mathfrak t}
\newcommand\p{\mathfrak p}

\newcommand\Sk{\operatorname{Sk}}

\newcommand\Spec{\operatorname{Spec}}

\newcommand\F{\operatorname{F}}
\newcommand\W{{\bf A}}
\newcommand\K{\mathbb K}
\newcommand\U{\mathcal U}

\newcommand\D{\mathcal D}
\newcommand\Ann{\operatorname{Ann}}
\newcommand\Id{\mathfrak{Id}}

\newcommand\Walg{\mathcal W}
\newcommand\Z{\mathbb Z}

\newcommand\M{\mathcal M}

\newcommand\N{\mathbb N}
\newcommand\gr{\operatorname{gr}}
\newcommand\Orb{\mathbb O}
\newcommand\Ocat{\mathcal O}
\newcommand\I{\mathcal I}
\newcommand\J{\mathcal J}
\renewcommand\sl{\mathfrak{sl}}

\newcommand\Sp{\mathop{\rm Sp}\nolimits}

\newcommand\gl{\mathfrak{gl}}
\newcommand\Span{\operatorname{Span}}
\newcommand\Hom{\operatorname{Hom}}
\newcommand{\ad}{\mathop{\rm ad}\nolimits}

\newcommand\Centr{\mathcal Z}

\newcommand{\VA}{\operatorname{V}}
\newcommand{\End}{\operatorname{End}}
\newcommand{\Kfun}{\mathcal{K}}

\newcommand\KF{\operatorname{K}}

\newcommand\Irr{\operatorname{Irr}}
\newcommand\Prim{\operatorname{Prim}}

\newcommand\SL{\operatorname{SL}}
\newcommand\Der{\operatorname{Der}}
\newtheorem{Thm}{Theorem}[section]
\newtheorem{Prop}[Thm]{Proposition}
\newtheorem{Cor}[Thm]{Corollary}

\theoremstyle{definition}

\newtheorem{Conj}[Thm]{Conjecture}
\numberwithin{equation}{section}
\title{Finite W-algebras}
\author{Ivan Losev}
\thanks{Supported by the NSF grant DMS-0900907}
\thanks{Address: Massachusetts Institute of Technology, Department of Mathematics,
77 Massachusetts Avenue, Cambridge MA 02139, USA. E-mail: ivanlosev@math.mit.edu}
\thanks{{\it 2010 Mathematics Subject Classification}: 
Primary 16G99,17B35; Secondary 53D20,53D55}
\thanks{{\it Keywords}: W-algebra,  semisimple Lie algebra, nilpotent orbit, universal enveloping algebra, primitive ideal,
Whittaker module.}

\begin{document}

\begin{abstract} A finite W-algebra is an associative algebra constructed from a semisimple Lie algebra
and its nilpotent element. In this survey we review recent developments in the representation
theory of  W-algebras. We emphasize various interactions between W-algebras and universal
enveloping algebras.
\end{abstract}

% Your AMS 200 Classification should come here %

% Any keywords %

% Do not remove next line %
\maketitle
\tableofcontents

% Now the text of your article starts
\section{Introduction}\label{SECTION_Intro}
Our base field $\K$ is supposed to be algebraically closed and of characteristic zero.

A finite W-algebra is an associative algebra constructed from a pair $(\g,e)$, where
$\g$ is a finite dimensional semisimple Lie algebra, and $e$ is a nilpotent element
of $\g$. A W-algebra should be thought as a generalization of the universal enveloping algebra
$U(\g)$. The latter can be considered as the W-algebra for the pair $(\g,0)$.

The study of W-algebras traces back to the celebrated paper \cite{Kostant} of Kostant.
This paper essentially treats the case when the element $e$ is principal (i.e., the adjoint orbit
of $e$ is dense in the nilpotent cone of $\g$). Kostant's motivation came basically from the study
of  Whittaker vectors and of Whittaker models. \cite{Kostant} was followed by the thesis
 \cite{Lynch} of Lynch who was a student of Kostant. In \cite{Lynch} Kostant's results were
(partially) generalized to arbitrary {\it even} nilpotent elements. During the 80's  Whittaker models (in the sense different from Kostant's) were also considered in %\cite{Kawanaka1},\cite{Kawanaka2} (in the context of finite Chevalley groups),
\cite{Moeglin1},\cite{Moeglin2},
where they were applied to the  study of certain primitive ideals in $U(\g)$. %(in the context of universal enveloping algebras).

In the 90's finite W-algebras attracted some attention from mathematical physicists,
see, for example, \cite{BT},\cite{RS},\cite{VD}. One of the main motivations
for their interest was a relationship between finite and {\it affine} W-algebras. The latter
are certain vertex algebras modeling the so called W-symmetry from Conformal field
theory.

In \cite{Premet1} Premet gave a general definition of a W-algebra. Premet's interest to
the subject was motivated by the study of non-restricted representations of semisimple
Lie algebras in positive characteristic. The  paper \cite{Premet1} initiated a lot of work
on different, mostly representation theoretic, aspects of W-algebras.

Apart from being of independent  interest, finite W-algebras have several connections to other objects studied in Representation theory. Let us summarize these connections.

A) It seems that the most straightforward connection is to the universal enveloping algebras
of semisimple Lie algebras.   This connection can be informally explained as follows.
According to the Orbit method, to an infinite dimensional representation of $\g$ one should be able
to assign a nilpotent orbit in $\g^*(\cong \g)$. For instance, to a "nice" (e.g., irreducible)
Harish-Chandra $\g$-bimodule one assigns a dense orbit in its associated variety.
Then there is a hope (that sometimes converts into proofs) that one can reduce
the study of an infinite dimensional $\g$-module to the study of a {\it finite dimensional}
module over the W-algebra corresponding to the nilpotent orbit in interest.
A relationship between the W-algebras and $U(\g)$ is studied, for example, in \cite{BGK},\cite{Ginzburg_HC},\cite{Wquant}-\cite{Miura},\cite{Premet2}-\cite{Premet4},\cite{Skryabin}.

B) There is a connection between the representation theory of W-algebras (in characteristic zero)
and that of semisimple Lie algebras in positive characteristic. In a sentence,  any reduced
enveloping algebra turns out to be Morita equivalent to an appropriate reduced W-algebra.
One can relate the representation theories of W-algebras in positive and in zero characteristics.
This relationship was successfully used in Premet's papers, see \cite{Premet0},\cite{Premet1},\cite{Premet3},\cite{Premet4}.

C) For classical Lie algebras there is a connection between W-algebras and (twisted) Yangians.
This connection was first discovered in \cite{RS} and then studied further in
\cite{Brown1},\cite{Brown2},\cite{BK1},\cite{BK2},\cite{R}. Also W-algebras are related to the
(cyclotomic quotients of)  degenerate affine Hecke algebras, \cite{BK3}.

D) As we mentioned above, finite W-algebras are related to their affine counterparts.
This relation can be made formal. To any vertex algebra one can assign an associative algebra
called the {\it Zhu algebra}. The importance of the Zhu algebra is that its representation
theory controls much of the representation theory of the initial vertex algebra.
A finite W-algebra is closely related to the Zhu algebra of the corresponding affine W-algebra. For details
the reader is referred to \cite{DSK}.

In the present paper we are mostly interested in A). We  also briefly explain  C), while
B) remains almost untouched and we do not discuss D) at all. Therefore we suppress the adjective
"finite" while speaking about W-algebras.  Another review on
W-algebras \cite{Wang} by W. Wang have already appeared. Some topics not discussed
(or discussed very briefly) in our survey can be found there.

This paper is organized as follows. In Section \ref{SECTION_W_HR}
we discuss topics related to a definition of a W-algebra via Hamiltonian reduction,
which is essentially due to Premet. In Section \ref{SECTION_W_DQ}
we explain the definition of a W-algebra  based on the
Deformation quantization, \cite{Wquant}. The next three sections describe connections between
W-algebras and $U(\g)$. In Section \ref{SECTION_equivalences}
we discuss category equivalences between certain categories of modules for W-algebras
and for $U(\g)$. Section \ref{SECTION_ideals} describes a relationship
between the sets of two-sided ideals in the two algebras. This description leads
to a (partial) classification of irreducible finite dimensional representations
of W-algebras in terms of primitive ideals of $U(\g)$. Section \ref{SECTION_one_dim}
deals with one-dimensional modules over W-algebras. Finally in Section \ref{SECTION_type_A}
we explain the connection C) above mostly for $\g$ of type A.

In the beginning of each section its content is described in more detail.

{\bf Acknowledgements.} First of all, I would like to thank J. Brundan, V. Ginzburg, S. Goodwin, A. Kleshchev,
and A. Premet for numerous inspiring discussions on W-algebras.
 I also thank J. Brundan, V. Ginzburg,
and A. Premet for their remarks on a preliminary version of this text.

%In more detail, let $\kf$
%be an algebraically closed field of positive characteristic. Then to $\g$ there corresponds the
%{\it restricted} Lie algebra $\g_{\kf}$. There is a natural embedding of the symmetric algebra
%$S(\g_{\kf})$ into the center of $U(\g_{\kf})$ turning $U(\g_{\kf})$. The algebra $U(\g_{\kf})$
%is a finite $S(\g_{\kf})$. Therefore to $\eta\in \g_{\kf}^*$ one can assign the fiber $U_\eta(\g_{\kf})$ of
%$U(\g_{\kf})$ at $\eta$. This is a finite dimensional algebra (called the reduced
%universal enveloping algebra). The study of the representation theory of $U_\eta(\g_{\kf})$
%can be to some extent reduced to the case when $\eta$ is nilpotent. One can define  W-algebras
%and its reduced versions in the positive characteristic case. Then  the reduced version
%of the W-algebra turns out to be Morita equivalent to $U_\eta(\g_{\kf})$. In a sentence,
%the modular representation theory of W-algebras is the same

\subsection*{Notation and conventions}
Throughout the paper $G$ is a connected reductive group, $\g$ is its Lie algebra.
We choose a nilpotent element $e\in \g$ and pick $h,f\in \g$ forming an $\sl_2$-triple
with $e$, i.e., $[h,e]=2e, [h,f]=-2f, [e,f]=h$. Let $\Orb$ denote the $G$-orbit of $e$.
Also we fix a $G$-invariant non-degenerate symmetric form $(\cdot,\cdot)$ on $\g$. Using this form,
we identify $\g$ with $\g^*$.

We write $\U$ for the universal enveloping algebra $U(\g)$ of $\g$. By $\Centr$
we denote the center of $\U$. This is a polynomial algebra.

Let us also list some standard notation used below.

\setlongtables
\begin{longtable}{p{2.5cm}p{12.5cm}}
$A^{op}$& the opposite algebra of an algebra $A$.\\
$\Ann_A(M)$& the annihilator of an $A$-module $M$.\\
$\operatorname{End}_A(M)$& the  algebra of endomorphisms of an $A$-module $M$.\\
$\operatorname{End}(M)$&$:=\operatorname{End}_\K(M)$.\\
$\gr V$& the associated graded vector space of a filtered vector space $V$.\\
$H^\circ$& the unit component of an algebraic group $H$.\\
$\K[X]$& the algebra of regular functions on a variety $X$.\\
$\K[X]^\wedge_Y$& the algebra of functions on the completion of a variety $X$
along a subvariety $Y$.\\
$T^*X$& the cotangent bundle of a smooth variety $X$.\\
$\VA(\M)$& the associated variety of a finitely generated $\U$-module $\M$.\\
$\z(\h)$& the center of a Lie algebra $\h$.\\
$\z_\h(\mathfrak{f})$& the centralizer of $\mathfrak{f}$ in a Lie algebra $\h$.
\end{longtable}

\section{W-algebras via Hamiltonian reduction}\label{SECTION_W_HR}
In this section we discuss developments leading to and related to Premet's definition of a
W-algebra given in \cite{Premet1}. The first such development is, of course, Kostant's work,
\cite{Kostant}, where the case of a principal nilpotent element was treated. We describe
(very few of) Kostant's results in Subsection \ref{SUBSECTION_Kostant}. Then in
Subsection \ref{SUBSECTION_Uge_even} we mention a generalization of Kostant's constructions
to the case of an even nilpotent element due to Lynch, \cite{Lynch}. In Subsection \ref{SUBSECTION_Uge}
we provide Premet's definition in the form of a quantum Hamiltonian reduction.
In Subsection \ref{SUBSECTION_Slodowy} we show that the "quasiclassical limit" of
a W-algebra is the algebra of functions on a {\it Slodowy slice}  that is a transverse
slice to a nilpotent orbit in $\g$ introduced in \cite{Slodowy}.  In Subsection \ref{SUBSECTION_Uge_ramifications} we mention several ramifications of Premet's definition and in Subsection \ref{SUBSECTION_Uge_remarks}
discuss some properties of W-algebras that can be proved using this definition.

\subsection{Kostant's results: the case of a principal nilpotent element}\label{SUBSECTION_Kostant}
In this subsection we will explain (some of) Kostant's results, Section 2 of \cite{Kostant}.

Suppose the nilpotent element $e$ is  principal.
Set \begin{equation}\label{eq:defi}\g(i):=\{\xi\in\g| [h,\xi]=i\xi\},
\p:=\bigoplus_{i\geqslant 0}\g(i), \m:=\bigoplus_{i<0}\g(i),\chi:=(e,\cdot).\end{equation}
%Of course, $\p$ is a Borel subalgebra and $\m$ is the opposite maximal nilpotent subalgebra.

Let us describe $\p,\m,e,\chi$ in more conventional terms.
Let $\h\subset\g$ be a Cartan subalgebra of $\g$, $\Delta\subset\h^*$ the corresponding root system,
and $\Pi$ a system of simple roots in $\Delta$. Further, for $\alpha\in \Delta$ let $e_\alpha$ denote
a corresponding weight vector in $\g$. Finally, let $\rho^\vee$ denote half the sum of all positive coroots
(=the sum of all fundamental co-weights). Replacing $(e,h,f)$ with a $G$-conjugate triple, we may assume that
$h=2\rho^\vee$ and $e=\sum_{\alpha\in \Pi} e_{\alpha}$. So
$\p$ becomes the positive Borel subalgebra $\borel\subset\g$, $\m$ becomes the negative maximal nilpotent subalgebra
$\n_-$, while $\chi$ is a non-degenerate character of $\m$.

Define the shift $ \m_\chi:=\{\xi-\langle\chi,\xi\rangle, \xi\in \m\}$ of $\m$.
Then, thanks to the PBW theorem,
we get \begin{equation}\label{eq:even_decomp}\U=U(\p)\oplus \U\m_\chi.\end{equation}

Using this decomposition,
one can define an action of $\m$ on $U(\p)$. Namely, identify
$U(\p)$ with the quotient $\U/\U\m_\chi$ using (\ref{eq:even_decomp}). The adjoint action of $\m$ on $\U$ descends to $\U/\U\m_\chi$. Using the identification $U(\p)\cong \U/\U\m_\chi$, we get an $\m$-action on $U(\p)$.

By definition, a W-algebra $U(\g,e)$ is the invariant subalgebra $U(\p)^\m$.
In other words, $U(\g,e)$ is the {\it quantum Hamiltonian reduction}
$$(\U/\U\m_\chi)^{\ad\m}:=\{a+\U\m_\chi: [\xi,a]\in \U\m_\chi,\forall \xi\in\m\}.$$
The multiplication on the last space is defined by $(a+\U\m_\chi)(b+\U\m_\chi)=ab+\U\m_\chi$.

It turns out that $U(\g,e)$ is naturally isomorphic to the center $\Centr$ of $\U$. Namely, the inclusion
$\Centr\hookrightarrow \U$ gives rise to a natural map $\Centr\rightarrow \U/\U\m_\chi$. Its image clearly consists of $\m$-invariants. So we get a homomorphism
$\Centr\rightarrow U(\g,e)$. By Theorem 2.4.1 in \cite{Kostant}, this homomorphism is an isomorphism. In particular,
we get an embedding of $\Centr$ into $U(\p)$. This embedding is of importance in the quantization of Toda
systems, see \cite{Kostant_Toda}.

\subsection{Generalization: the case of even $e$}\label{SUBSECTION_Uge_even}
Recall that $e$ is called {\it even} if all
eigenvalues of $\ad h$ on $\g$ are even.
Define $\g(i),\p,\m,\chi$ by (\ref{eq:defi}).
It is clear that $\p$ is a parabolic subalgebra of $\g$ and $\m$
is the nilpotent radical of the opposite parabolic. In \cite{Lynch}
Lynch generalized Kostant's definition and introduced an algebra $U(\g,e):=U(\p)^\m=
(\U/\U\m_\chi)^{\ad \m}$.

There is an embedding $U(\g,e)\hookrightarrow
U(\g(0))$ sometimes called the {\it generalized Miura transform}. It is obtained by restricting the natural projection $U(\p)\twoheadrightarrow U(\g(0))$ to $U(\g,e)\subset U(\p)$. The restriction is injective by \cite{Lynch}, Corollary 2.3.2.

\subsection{Definition of $U(\g,e)$: the general case}\label{SUBSECTION_Uge}
Now let $e\in\g$ be an arbitrary (nonzero) nilpotent element.
Let the decomposition
$\g=\bigoplus_{i\in \Z}\g(i)$ and the element $\chi\in\g^*$ be given by (\ref{eq:defi}).
Following Premet, \cite{Premet1}, we still
can define a W-algebra $U(\g,e)$ as the quantum Hamiltonian reduction $(\U/\U\m_\chi)^{\ad \m}$
provided we can define a suitable analog of the subalgebra $\m\subset\g$
considered in the previous subsection.

A subalgebra $\m$ we need is constructed as follows.
Consider a skew-symmetric form $\omega_\chi$ on $\g$ given by $\omega_\chi(\xi,\eta)=\langle\chi,[\xi,\eta]\rangle$.
It follows easily from the representation theory of $\sl_2$ that  the restriction of $\omega_\chi$ to the subspace $\g(-1)$ is non-degenerate. Pick a lagrangian subspace $l\subset \g(-1)$ and set $\m:=l\oplus \bigoplus_{i\leqslant -2}\g(i)$.

It is clear that $\m$ is a subalgebra in $\g$ consisting of nilpotent elements.
Also since $\omega_\chi$ vanishes on $l$, we see that $\langle\chi, [\m,\m]\rangle=0$.
So $\chi$  is indeed a character of $\m$.
%Finally, we want to make the following remark whose importance
%will become clear later. Consider the subspace $V:=\im\ad(f)\subset \g$. The restriction of
%$\omega_\chi$ to $V$ is non-degenerate so $V$ becomes a symplectic vector space. It is easy
%to check that $\m$ is a lagrangian subspace in $V$.

We set $U(\g,e):=(\U/\U\m_\chi)^{\ad \m}$. The reader should notice that, a priory, this definition
is ambiguous: $\m$ and hence $U(\g,e)$ depend on the choice of $l$. However, we will see in Subsection
\ref{SUBSECTION_Uge_ramifications} that two $W$-algebras constructed using different choices of $l$
are canonically isomorphic.

We finish the subsection with a few historical remarks.
Direct analogs of $\m$ and of its shift $\m_\chi$ in the setting of finite Chevalley groups first appeared in \cite{Kawanaka1}. Then Moeglin used $\m_\chi$ to define "Whittaker models" for primitive
ideals of $U(\g)$, \cite{Moeglin1},\cite{Moeglin2}. We will describe her results in more detail in Section \ref{SECTION_one_dim}.
Later the subalgebra $\m$ played an important role in Premet's proof of the Kac-Weisfeller conjecture
on the dimension of a non-restricted representation of a semisimple Lie algebra in positive characteristic,
see \cite{Premet0}.

\subsection{Classical counterpart: the Slodowy slice}\label{SUBSECTION_Slodowy}
The algebra $U(\g,e)$ has an interesting filtration, called the {\it Kazhdan}
filtration.

To define it we first introduce a new filtration on $\U$. Recall that the algebra $\U$ has the standard,
PBW filtration: the subspace  $\F^{st}_i\U$ of elements of degree $\leqslant i$, by definition, is spanned
by all monomials $\xi_1\ldots\xi_j, j\leqslant i, \xi_1,\ldots,\xi_j\in\g$. For $j\in \Z$ set
$\U(j):=\{u\in \U| [h,u]=ju\}$. Define the {\it Kazhdan filtration} $\KF_i \U$ on $\U$ by
$\KF_i\U:=\sum_{2j+k\leqslant i}\F^{st}_j\U\cap \U(k)$. We remark that the associated graded algebra
of $\U$ with respect to the Kazhdan filtration is still naturally isomorphic to the symmetric algebra
$S(\g)$.

Being a subquotient of $\U$, the algebra $U(\g,e)$ has a Kazhdan filtration $\KF_i U(\g,e)$ inherited
from $\U$. We remark that $\KF_{0}\U\subset \K+\U\m_\chi$ so the Kazhdan filtration on $U(\g,e)$ is positive
in the sense that $\KF_0 U(\g,e)=\K$.

 It turns out that
 the associated graded algebra $\gr U(\g,e)$ of $U(\g,e)$  is naturally isomorphic to  the algebra of functions on the   {\it Slodowy slice} $S:=e+\ker\ad(f)$, \cite{Slodowy} (in the case when $e$ is principal $S$ appeared in
\cite{Kostant_pol}).  It follows from the representation theory of $\sl_2$ that $S$ is transverse to  $\Orb$. In the sequel it will be convenient for us to consider $S$ as an affine subspace in $\g^*$ via the identification $\g\cong\g^*$. In particular, $\chi\in S$.

We need an action of the one-dimensional torus $\K^\times$ on $\g^*$ that stabilizes $S$ and contracts it
to $\chi$. Namely, the $\sl_2$-triple $(e,h,f)$ defines a homomorphism $\SL_2(\K)\rightarrow G$. The group
$\K^\times$ is embedded into $\SL_2(\K)$ via $t\mapsto \operatorname{diag}(t,t^{-1})$. Composing these two
homomorphisms we get a homomorphism (in fact, an embedding) $\gamma:\K^\times\rightarrow G$. For
$\xi\in \g(i)$ we have $\gamma(t).\xi= t^i\xi$. Consider a $\K^\times$-action (called the Kazhdan action) on $\g^*$ given by
$t\cdot\alpha= t^{-2}\gamma(t) \alpha$. This action fixes $\chi$. Also it is easy to see that
it preserves $S$. Finally, the representation theory of $\sl_2$ implies that $\ker\ad(f)\subset \bigoplus_{i\leqslant 0}\g(i)$. It follows that the action of $\K^\times$ contracts $S$ to $\chi$: $\lim_{t\rightarrow \infty} t.s=\chi$
for any $s\in S$. The contraction property has several very nice corollaries. For example, $S$ intersects
an adjoint orbit $\Orb'$ if and only if $\Orb\subset \overline{\K^\times \Orb'}$ and in this case
the intersection $S\cap\Orb'$ is transversal.

The Kazhdan action gives rise to a (positive) grading on  the algebra $\K[S]$ of regular functions on $S$.
The following result  was essentially obtained by Premet, \cite{Premet1}, Theorem 4.6
(Kostant and Lynch also proved this in the special cases they considered).

\begin{Thm}\label{Thm:ass_gr}
There is an isomorphism $\gr U(\g,e)\cong \K[S]$ of graded algebras.
\end{Thm}

As was shown by Gan and Ginzburg, \cite{GG}, this result is a manifestation of the "quantization commutes with reduction" principle, see Subsection \ref{SUBSECTION_Uge_ramifications} for details.

%We remark that in the special case, where $e$ is principal, all the constructions
%and results of this subsection are due    to Kostant, \cite{Kostant}.
%They were generalized to an arbitrary even element  by Lynch. In fact, the approach in \cite{GG}
%basically follows Kostant's.

\subsection{Ramifications}\label{SUBSECTION_Uge_ramifications}
%First of all, let us make a trivial remark. Although above we considered a semisimple Lie algebra
%$\g$ everything works just as well with a reductive Lie algebra. Clearly, $U(\g,e)=S(\z(\g))\otimes U([\g,\g])$
%
First, let us mention the work of Gan and Ginzburg, \cite{GG}, where they gave a ramification
of Premet's definition showing, in particular, that $U(\g,e)$ does not depend on the choice of $l\subset \g(-1)$.

Namely, let $l\subset \g(-1)$ be an arbitrary isotropic subspace of $\g(-1)$ (e.g., $\{0\}$).
Let $l^\angle$ denote the skew-orthogonal complement to $l$ in $\g(-1)$. Set $\m^l:=l\oplus \bigoplus_{i\leqslant -2}
\g(i), \n^l:=l^\angle\oplus \bigoplus_{i\leqslant -2}\g(i)$. Then $\m^l\subset \n^l$, $\n^l$ consists
of nilpotent elements, and $\langle\chi, [\m_l,\n_l]\rangle=0$. Let $N^l$ be the connected subgroup of
$G$ with Lie algebra $\n^l$. Then $\m^l$ and the character $\chi:\m^l\rightarrow \K$ are stable under the
adjoint action of $N^l$. So $N^l$ acts naturally on the quotient $\U/\U\m^l_\chi$, where
$\m^l_\chi:=\{\xi-\langle\chi,\xi\rangle, \xi\in \m^l\}$.
Let $U(\g,e)^l:=(\U/\U\m^l_\chi)^{N^l}$ be the space
of invariants. It is easy to check that it  has a natural algebra
structure. It also has a Kazhdan filtration $\KF_i U(\g,e)^l$, compare with the previous subsection.

Now let us remark that for $l_1\subset l_2$ we have a natural $\U$-module homomorphism
$\U/\U\m^{l_1}_\chi\rightarrow \U/\U \m^{l_2}_\chi$ that gives rise to a filtered algebra
homomorphism $U(\g,e)^{l_1}\rightarrow U(\g,e)^{l_2}$. It turns out that the latter is an isomorphism.

Also Gan and Ginzburg gave a very transparent explanation of an isomorphism $\gr U(\g,e)\cong \K[S]$.
Namely, consider the restriction map $\pi:\g^*\rightarrow \m^{l*}$. The affine subspace $\pi^{-1}(\chi|_{\m^l})\subset \g^*$ is $N^l$-stable. Also it is easy to see that $S\subset \pi^{-1}(\chi|_{\m^l})$. So we can consider a morphism
$N^l\times S\rightarrow \pi^{-1}(\chi|_{\m^l}), (n,s)\mapsto ns$.
According to \cite{GG}, this is an isomorphism (of algebraic
varieties). Therefore $\K[S]$ gets identified with the classical Hamiltonian
reduction  $(S(\g)/S(\g)\m^l_\chi)^{N^l}$.

Another ramification of the original definition of $U(\g,e)$ comes from the notion of a
good grading on $\g$,  \cite{EK}. A grading $\g=\bigoplus_{i\in \Z}\g(i)$ is said to be {\it good}
for  $e$ if $e\in \g(2)$ and  $\ker\ad(e)\subset \bigoplus_{i\geqslant 0}\g(i)$.
For instance, the grading given by (\ref{eq:defi})  is good. For a
comprehensive study of good gradings see \cite{EK}.

Given a good grading on $\g$, one constructs $\m\subset\g$ and defines
$U(\g,e)$ using  $\m$ analogously to the above. The algebra $U(\g,e)$ does
not depend on the choice of a good grading. This was first proved in \cite{BG}.

The definition involving an arbitrary good grading is often useful. For example, one can sometimes
find an {\it even} good grading when  $e$ is not even itself and embed $U(\g,e)$ into $U(\p)$
for an appropriate parabolic subalgebra $\p\subset\g$,
compare with Subsection \ref{SUBSECTION_Uge_even}. This is always the case when
$\g\cong \sl_n$, see \cite{BK1}, Introduction.

Also it is worth mentioning that there is a related definition of $U(\g,e)$ via the BRST quantization procedure
which was used by physicists in the 90-s, see \cite{BT}. The proof that the BRST definition is equivalent to
the one given above was obtained in \cite{DSK_appendix}. See also \cite{Wang}, Section 3.

\subsection{Additional properties of $U(\g,e)$}\label{SUBSECTION_Uge_remarks}
We want to make a few remarks about other properties of $U(\g,e)$.

Recall that $\Centr$ stands for the center of $\U$.
Restricting the natural map
$\U^{\ad\m}\rightarrow U(\g,e)$ to $\Centr\subset \U^{\ad\m}$, we get an algebra homomorphism $\Centr\rightarrow U(\g,e)$.
By \cite{Premet1}, 6.2, this homomorphism is an embedding. It is clear that the image of
$\Centr$ lies in the center of $U(\g,e)$. Further, according to the footnote to
Question 5.1 in \cite{Premet2}, the image of $\Centr$ actually coincides with
the center of $U(\g,e)$ (Premet attributes the proof to Ginzburg). This generalizes Kostant's result mentioned in Subsection \ref{SUBSECTION_Kostant}.

Also we remark that there is a natural action of the group $Q:=Z_G(e,h,f)$ on $U(\g,e)$. Namely, take $l=\{0\}$ in the Gan and Ginzburg definition. Then $Q$ stabilizes both $\m^l_\chi$ and $N^l$ and so
acts on $U(\g,e)^l$. Let $\q$ stand for the Lie algebra of $Q$. In \cite{Premet2} Premet constructed a Lie algebra embedding  $\q\hookrightarrow U(\g,e)$ such that the adjoint action of $\q$ on $U(\g,e)$ coincides with the
differential of the $Q$-action.

\section{W-algebras via Deformation quantization}\label{SECTION_W_DQ}
In this section we review the definition of W-algebras from \cite{Wquant}. It
is based on Deformation quantization: a W-algebra is realized as an algebra of $G$-invariants
in a quantization of a certain affine symplectic $G$-variety (called an {\it equivariant Slodowy
slice}). In Subsection \ref{SUBSECTION_Fedosov} we briefly explain generalities on
star-products and on Fedosov's method to construct them. In Subsection \ref{SUBSECTION_equiv_Slodowy}
we present  constructions of  equivariant Slodowy slices and of  W-algebras. Finally,
in Subsection \ref{SUBSECTION_Decomp} we present a very important basic result on W-algebras,
 the {\it decomposition theorem}.

\subsection{Fedosov quantization}\label{SUBSECTION_Fedosov}
In this subsection $X$ is a smooth affine variety equipped with a symplectic form $\omega$.
Let $\{\cdot,\cdot\}$ denote  the Poisson bracket on $\K[X]$
induced by $\omega$.
Let  a reductive group $\widetilde{G}$ act on $X$ preserving $\omega$.
By $\xi_X$ we denote  the image of $\xi\in\widetilde{\g}$ under the homomorphism $\widetilde{\g}\rightarrow \Der(\K[X])$
induced by the action.

We suppose that the $\widetilde{G}$-action is Hamiltonian, that is, admits a {\it moment map} $\mu:X\rightarrow \widetilde{\g}^*$, i.e.,  a $\widetilde{G}$-equivariant morphism having the following property: for $H_\xi:=\mu^*(\xi),\xi\in\widetilde{\g},$
we have $\{H_\xi,\cdot\}=\xi_X$. Finally, we suppose that $X$ is equipped with a $\K^\times$-action that commutes
with $\widetilde{G}$ and satisfies $t.\omega= t^2\omega, t.H_\xi=t^2 H_\xi$ for all $t\in \K^\times,\xi\in\widetilde{\g}$. We will
present  examples of this situation below.

By a {\it star-product} on $\K[X]$ (or on $X$) we mean a
$\K$-bilinear map $*:\K[X]\times \K[X]\rightarrow \K[X][[\hbar]],
(f,g)\mapsto f*g:=\sum_{i=0}^\infty D_i(f,g)\hbar^{2i}$ satisfying the following axioms:
\begin{itemize}
\item[(a)] The associativity axiom: a natural  extension of
$*$ to a $\K[[\hbar]]$-bilinear map
$\K[X][[\hbar]]\times \K[X][[\hbar]]\rightarrow \K[X][[\hbar]]$ is an associative product,
and $1\in \K[X]\subset \K[X][[\hbar]]$ is a unit for $*$.
\item[(b)] The compatibility axiom: $D_0(f,g)=fg, D_1(f,g)-D_1(g,f)=\{f,g\}$. Equivalently,
$f*g\equiv fg \mod \hbar^2$ and $[f,g]\equiv \hbar^2 \{f,g\}\mod \hbar^4$.
\item[(c)] The locality axiom: $D_i$ is a bidifferential operator of order at most $i$
(i.e., for any fixed $f$ the map $\K[X]\rightarrow \K[X]: g\mapsto D_i(f,g),$ is a differential operator
of order at most $i$, and the same for any fixed $g$).
\end{itemize}
When we consider $\K[X][[\hbar]]$ as an algebra with respect to the star-product, we
call it a {\it quantum algebra}.

We remark that the usual definition of a star-product looks like $f*g=\sum_{i=0}^\infty D_i(f,g)\hbar^i$
and in our definition we have $\hbar^2$ instead of $\hbar$. The reason for this ramification
is that our version is better compatible with the {\it Rees algebra construction}. This construction
allows to pass from filtered $\K$-algebras to graded $\K[\hbar]$-algebras.

%When we consider $\K[X][[\hbar]]$ as an algebra with respect to the star-product, we call it a quantum
%algebra.

We will also need  $*$ to be compatible with the $\widetilde{G}$- and $\K^\times$- actions on $X$.
\begin{itemize}
\item[(d)] $\widetilde{G}$-invariance: $D_i:\K[X]\otimes \K[X]\rightarrow \K[X]$ is $\widetilde{G}$-equivariant.
\item[(e)] Homogeneity: $D_i$ has degree $-2i$ with respect to $\K^\times$: i.e., for
$f,g\in \K[X]$ of degrees $j,k$  the element
$D_i(f,g)$ has degree $k+j-2i$.
\end{itemize}

Under the conditions (d) and (e), the product $\widetilde{G}\times \K^\times$ acts on $\K[X][[\hbar]]$ by automorphisms
with $g.\hbar=\hbar, t.\hbar=t\hbar$ for all $g\in G, t\in \K^\times$.

It turns out that a star-product on $X$ satisfying additionally (d) and  (e) always exists.
It is provided, for example, by Fedosov's construction, \cite{Fedosov1},\cite{Fedosov2}.
Fedosov constructed a star-product on a $C^\infty$-manifold starting from a {\it symplectic connection} $\nabla$ and a closed $\K[[\hbar^2]]$-valued form $\Omega$.
By definition, a symplectic connection is  a torsion-free connection on the tangent bundle such that the symplectic form is flat. Fedosov's construction can be carried over to the algebraic
setting as long as a variety in consideration admits a symplectic connection. Since $X$ is affine,
this is the case, and, moreover, one can, in addition, assume that a symplectic connection
is $\widetilde{G}\times \K^\times$-invariant, see \cite{Wquant}, Proposition 2.2.2. For our purposes,
it will be enough to consider the original construction from \cite{Fedosov1}, where
$\Omega$ is not used (i.e., equals 0).

The following proposition follows from results of Fedosov, see \cite{HC}, Theorem 2.1.2 for details.

\begin{Prop}\label{Prop:quantization}
Let $X$ be as above, and $\nabla$ be a $\widetilde{G}\times \K^\times$-invariant symplectic connection on $X$.
Further, let $*$ be the star-product produced from $\nabla$ by the Fedosov construction. Then
$*$ is $\widetilde{G}$-invariant and homogeneous. Moreover, the map $\xi\mapsto H_\xi$ is a {\rm quantum comoment
map} for the $\widetilde{G}$-action on $\K[X][[\hbar]]$, i.e., $\frac{1}{\hbar^2}[H_\xi, f]= \xi_X f$ for all
$f\in \K[X][[\hbar]],\xi\in\widetilde{\g}$.
\end{Prop}

Also, according to Fedosov, $*$ does not depend on the choice of $\nabla$ up to a suitably
understood isomorphism, see, for example, \cite{Wquant}, Proposition 2.2.5, for details.

Let us consider two standard examples.

The first example is easy. Let $V$ be a vector space equipped with a non-degenerate form $\omega\in \bigwedge^2 V^*$.
Let $\widetilde{G}$ act on $V$ via a homomorphism $\widetilde{G}\rightarrow \Sp(V)$. Pick a homomorphism $\beta: \K^\times\rightarrow \Sp(V)^{\widetilde{G}}$ and define a $\K^\times$-action on $V^*$ by $t.\alpha=t^{-1}\beta(t)\alpha$. So we get a symplectic variety $X=V^*$ equipped with a $\widetilde{G}\times\K^\times$-action satisfying the assumptions above with the moment map
given by $\langle\mu(v),\xi\rangle=\frac{1}{2}\omega(\xi v,v)$. The algebra $\K[V^*]$ has a standard star-product called the {\it Moyal-Weyl} product. Namely, for $f,g\in \K[V^*]$ set $f*g:=m(\exp(\frac{\omega}{2}\hbar^2)f\otimes g)$. Here $m:\K[V^*]\otimes\K[V^*]\rightarrow \K[V^*]$ stands for the multiplication map, while $\omega\in \bigwedge^2 V^*$ is assumed to act on
$\K[V^*]\otimes \K[V^*]$  via contraction. The quantum algebra $\K[V^*][\hbar]$ is naturally identified with
the "homogeneous" version $\W_\hbar$ of the Weyl algebra of $V$, $\W_\hbar:=T(V)[\hbar]/(u\otimes v-v\otimes u-\hbar^2\omega(u,v), u,v\in V)$.

Our second example is more involved although is also standard.

Let $G$ be a connected reductive algebraic group.
 The cotangent bundle $X:=T^*G$ of $G$
is equipped with a canonical symplectic form $\omega$. Set $\widetilde{G}:=G\times G$ and consider the $\widetilde{G}$-action on $X$ induced from the two-sided action of $\widetilde{G}$ on $G$. In more detail, we can identify
$T^*G$ with $G\times\g^*$ using the trivialization by left-invariant forms.
Then the "left" action of $G$ on $X$ is given by $g.(g_1,\alpha)=(gg_1,\alpha)$,
while the "right" action is $g.(g_1,\alpha)=(g_1g^{-1}, g.\alpha)$.
Finally, let $\K^\times$ act on $X$  by $t.(g_1,\alpha)=(g_1,t^{-2}\alpha)$.
Clearly, $\omega$ is
$\widetilde{G}$-invariant and $t.\omega=t^2\omega$. A moment map
$\mu:X\rightarrow \widetilde{\g}^*=\g^*\oplus \g^*$ is given by $(g,\alpha)\mapsto
(g.\alpha,\alpha)$.

Pick a $\widetilde{G}\times \K^\times$-invariant connection $\nabla$ on $X$ and produce the star-product
$*$ from $\nabla$. From the grading considerations, we see that $\K[X][\hbar]$ is a subalgebra in
the quantum algebra $\K[X][[\hbar]]$.

There  is a standard alternative description of $\K[X][\hbar]$, see, for example, Subsection 7.1 of \cite{Miura}. Consider the algebra $\D(G)$ of linear differential operators on $G$.
Let  $\F_i\D(G)$  be the space of differential operators of order $\leqslant i/2$. Consider the Rees algebra $\D_\hbar(G):=\bigoplus_{i=0}^\infty  \F_i\D(G)\hbar^i\subset\D(G)[\hbar]$ of $\D(G)$.
Then there is a $\widetilde{G}\times \K^\times$-equivariant isomorphism $\K[X][\hbar]\cong \D_\hbar(G)$
of $\K[\hbar]$-algebras.

Taking the  $G$-invariants in
the algebra  $\K[T^*G][\hbar]\cong \D_\hbar(G)$
(say for the left $G$-action), we get a new (star-)product on $\K[\g^*][\hbar]=\K[T^*G][\hbar]^G$. But $\D_\hbar(G)^G$ is nothing else but a homogeneous version $\U_\hbar$ of the universal enveloping algebra $\U$ of $\g$, $\U_\hbar:=T(\g)[\hbar]/(\xi\otimes\eta-\eta\otimes\xi-\hbar^2[\xi,\eta], \xi,\eta\in \g)$.
In the next subsection we will use a similar recipe to define a W-algebra.

\subsection{Equivariant Slodowy slices and W-algebras}\label{SUBSECTION_equiv_Slodowy}
A variety we need in the approach to  W-algebras from \cite{Wquant} is as follows.
Recall  the Slodowy slice $S\subset \g^*$, Subsection \ref{SUBSECTION_Slodowy}.
Set $X:=G\times S\subset G\times \g^*=T^*G$. The variety $X$ is called the {\it equivariant
Slodowy slice}. Clearly, $X\subset T^*G$ is stable with respect to the left $G$-action. Also it is stable under the
restriction of the right $G$-action to $Q=Z_G(e,h,f)$. Finally, $X$ is stable under a {\it Kazhdan}
$\K^\times$-action given by $t.(g,\alpha)=(g\gamma(t)^{-1}, t^{-2}\gamma(t)\alpha)$, where
$\gamma:\K^\times \rightarrow G $ was introduced in Subsection \ref{SUBSECTION_Slodowy}.
Consider the 2-form $\omega$ on $X$ obtained by the restriction of the natural symplectic
form from $T^*G$. One can show that $\omega$ is non-degenerate.
So $X$ becomes a symplectic variety. It satisfies the assumptions in the beginning of the previous subsection
with $\widetilde{G}:=G\times Q$, the Kazhdan action of $\K^\times$ and a moment map $X\rightarrow \g^*\oplus \q^*$
restricted from $T^*G$.

Pick a $\widetilde{G}\times \K^\times$-invariant symplectic connection
$\nabla$ on $X$ and produce a star-product $f*g=\sum_{i=0}^\infty D_i(f,g)\hbar^{2i}$,
using the Fedosov construction.  \cite{Wquant}, Proposition 2.1.5 implies
 that $\K[X][\hbar]\subset \K[X][[\hbar]]$ is closed with respect to the star-product.  We call the  quantum algebra $\K[X][\hbar]$  a {\it homogeneous equivariant W-algebra} and denote it by $\widetilde{\Walg}_\hbar$.
A {\it homogeneous W-algebra} is, by definition, $\Walg_\hbar:=\widetilde{\Walg}_\hbar^G$. Finally,
define a W-algebra $\Walg$ as $\Walg_\hbar/(\hbar-1)$. So, as a vector space $\Walg$ is the same as  $\K[S]$
but the product on $\Walg$ is given by $fg:=\sum_{i=0}^\infty D_i(f,g)$.

The algebra $\Walg$ comes equipped with
\begin{itemize}
\item a filtration  $\F_i\Walg$ induced from the grading on $\Walg_\hbar$.
\item an action of $Q$.
\item a homomorphism (in fact an embedding) $\q\hookrightarrow \Walg$ of Lie algebras
such that the adjoint action of $\q$ on $\Walg$ coincides with the differential
of the $Q$-action.
\item a homomorphism $\Centr\rightarrow \Walg$ (induced from the quantum comoment
map $\g\rightarrow \widetilde{\Walg}_\hbar$).
\end{itemize}

It turns out that $\Walg$ is isomorphic to $U(\g,e)$. More precisely, we have the following result.

\begin{Thm}[\cite{Wquant}, Corollary 3.3.3]\label{Thm:iso}
There is a  filtration preserving  isomorphism $\Walg\rightarrow U(\g,e)$.
\end{Thm}

One can prove, in addition,
that this isomorphism is $Q$-equivariant (although this is not written down explicitly)
and intertwines the homomorphisms $\Centr\rightarrow \Walg, U(\g,e)$ (this is proved in \cite{HC},
the end of Subsection 2.2).

%Let us sketch one of the proofs of Theorem \ref{Thm:iso}. For this we need a construction
%of $\widetilde{\Walg}$ along the lines of Section \ref{SECTION_W_HR}. Consider the algebra
%$\D(G)$ of differential operators on $G$. We have two embeddings $\U\rightarrow \D(G)$
%induced by the left (denoted by $\iota_l$) and the right ($\iota_r$)
%actions of $G$ on itself. Set $\D(G,e):=(\D(G)/\D(G)\iota_r(\m_\chi))^{\ad \iota_r(\m)}$.
%This algebra comes   equipped with a $G$-action, a homomorphism $\U\rightarrow \D(G,e)$.

\subsection{Decomposition theorem}\label{SUBSECTION_Decomp}
%Note that usually a star-product is written as $f*g=\sum_{i=0}^\infty D_i(f,g)\hbar^i$ with
%$f*g-g*f-\hbar\{f,g\}\in\hbar^2 A[[\hbar]]$. The reason why we use $\hbar^2$ instead of
%$\hbar$ is that this choice is compatible with the Rees algebra construction, which is used
%to pass from a filtered $\K$-algebra to a graded $\K[\hbar]$-algebra.

Let $x$ denote the point $(1,\chi)\in X\subset T^*G=G\times\g^*$. We remark that the orbit $Gx$ is closed
(as any orbit in $T^*G$) and also $Q\times \K^\times$-stable.  %The star-products on $T^*G, X,V^*$ are written by means %of bidifferential operators.
Consider the formal neighborhoods $(T^*G)^\wedge_{Gx},X^\wedge_{Gx}$ of $Gx$ in $T^*G$ and $X$ and  the formal neighborhood $(V^*)^{\wedge}_0$ of $0$ in $V^*$.
Being defined by bidifferential operators, the star-products on $\K[T^*G][\hbar], \K[X][\hbar],\K[V^*][\hbar]$ extend to the corresponding completions $\K[T^*G]^\wedge_{Gx}[[\hbar]], \K[X]^\wedge_{Gx}[[\hbar]], $ $ \W_\hbar^\wedge:=\K[V^*]^\wedge_0[[\hbar]]$.

Taking the $G$-invariants in $\K[T^*G]^\wedge_{Gx}[[\hbar]],\K[X]^\wedge_{Gx}[[\hbar]]$ we get star-products
on the completions $\U^\wedge_\hbar:=\K[\g^*]^\wedge_\chi, \Walg^\wedge_\hbar:=\K[S]^\wedge_\chi[[\hbar]]$.
The algebras $\U^\wedge_\hbar,\Walg^\wedge_\hbar,\W^\wedge_\hbar$ come equipped with natural (complete and
separated) topologies.
We remark that the completions $\U^\wedge_\hbar, \Walg^\wedge_\hbar,\W^\wedge_\hbar$ can be defined
completely algebraically, as the inverse limits of $\U_\hbar,\Walg_\hbar, \W_\hbar$ with respect to
the powers of appropriate maximal ideals, see \cite{HC}, Subsection 2.4 for details.

The following theorem follows from \cite{Wquant}, Theorem 3.3.1.

\begin{Thm}\label{Thm:decomposition}
There is a $Q\times \K^\times$-equivariant  isomorphism $\Phi_\hbar:\U^\wedge_\hbar\rightarrow
\W^\wedge_\hbar\widehat{\otimes}_{\K[[\hbar]]}\Walg^\wedge_\hbar$ of topological $\K[[\hbar]]$-algebras.
\end{Thm}

Here $\widehat{\otimes}$ stands for the completed tensor product: we take the usual tensor product
of topological $\K[[\hbar]]$-algebras and then complete it with respect to the induced topology.

Theorem \ref{Thm:decomposition} is extremely important in the study of $W$-algebras.
It can be used to prove Theorem \ref{Thm:iso}, to prove the category equivalence theorems
\ref{Thm:Skryabin},\ref{Thm_Ocat} in the next section,  and also to relate the sets of two-sided ideals of
$\U$ and of $\Walg$, see Section \ref{SECTION_ideals}.

%\subsection{Equivalence of definitions}\label{SUBSECTION_equivalence}

\section{Category equivalences}\label{SECTION_equivalences}
This section is devoted to the description of two category equivalences
between suitable categories of $\Walg$-modules and of $\U$-modules.
In the first subsection we recall an equivalence proved by Skryabin
in \cite{Skryabin}. This is an equivalence between the category of all
$\Walg$-modules and the category of {\it Whittaker} $\U$-modules.
Then we discuss some corollaries of Skryabin's theorem, in particular, a localization
theorem due to Ginzburg, \cite{Ginzburg_HC}.
Subsection \ref{SUBSECTION_Ocat} deals with a ramification of Skryabin's equivalence
conjectured in \cite{BGK} and proved in \cite{LOCat}. This is an equivalence between
the {\it category $\Ocat$} for a W-algebra  and the category of {\it generalized
Whittaker} $\U$-modules.

\subsection{Whittaker modules and Skryabin's equivalence}\label{SUBSECTION_Skryabin}
Recall that in Subsection \ref{SUBSECTION_Uge} we have defined the W-algebra
$\Walg=U(\g,e)$ as the quantum Hamiltonian reduction $(\U/\U\m_\chi)^{\ad\m}$.
In other words, $\Walg=\End_\U(\U/\U\m_\chi)^{op}$. In particular, $\U/\U\m_\chi$ is a $\U$-$\Walg$-bimodule.
%The action of $\m_\chi$ by left multiplication  on $\U/\U\m_\chi$ coincides with the action
%induced from the adjoint $\m_\chi$-action on $\U$ and hence
%is locally nilpotent.

We say that a $\U$-module $M$ is {\it Whittaker} if the action of $\m_\chi$ on $M$ is locally nilpotent.
For instance, $\U/\U\m_\chi$ is easily seen to be Whittaker.
Whittaker modules form a Serre subcategory  in the category
$\U$-$\operatorname{Mod}$. Denote the category of Whittaker $\U$-modules by $\operatorname{Wh}$.

The bimodule $\U/\U\m_\chi$ gives rise to the following functors:
\begin{align*}
\operatorname{Wh}\rightarrow \Walg\text{-}\operatorname{Mod}&: M\mapsto \Hom_\U(\U/\U\m_\chi,M)=M^{\m_\chi}:=\\ &:=\{m\in M: \xi m=\langle\chi,\xi\rangle m, \forall \xi\in\m\}.\\
\Walg\text{-}\operatorname{Mod}\rightarrow \operatorname{Wh}&:
N\mapsto \U/\U\m_\chi\otimes_{\Walg} N.
\end{align*}
We denote the second functor by $\Sk$.

The following important theorem was proved in \cite{Skryabin}.
\begin{Thm}\label{Thm:Skryabin}
The functors above are quasi-inverse equivalences.
\end{Thm}

Let us mention several important corollaries of this theorem.

The Beilinson-Bernstein localization
theorem, \cite{BB}, is a crucial result in the representation theory of $\U$. There is an analog of this theorem for W-algebras due to Ginzburg, \cite{Ginzburg_HC}.
See also \cite{DK} for an alternative approach.

Recall the Beilinson-Bernstein theorem. Pick a Cartan subalgebra $\h\subset\g$. Let $\Delta\subset\h^*$
be the root system, $W$ the Weyl group, and $\Pi\subset \Delta$ be a system of simple roots.
Recall the dot action of $W$ on $\h^*$ given by $w\cdot\lambda= w(\lambda+\rho)-\rho$, where,
as usual, $\rho$ stands for the half of the sum of all positive roots. The center $\Centr$ of $\U$
gets identified via the Harish-Chandra isomorphism with the invariant algebra $\K[\h^*]^W$,
the invariants are taken with respect to the dot action.

To any $\lambda\in\h^*$ one assigns a sheaf $\Dcal_\lambda$ of twisted differential operators on
the flag variety $\Fl$ of $G$, see  \cite{BB}. The algebra $\Gamma(\Fl,\Dcal_\lambda)$ of global
sections is naturally identified with the quotient $\U_{\lambda}:=\U/\U I_\lambda$, where $I_\lambda$
 denotes the maximal ideal of $W\cdot \lambda$ in $\Centr$.
So to a $\Dcal_\lambda$-module $M$ one can assign the $\U_\lambda$-module
$\Gamma(\Fl,M)$. The functor $\Gamma(\Fl,\bullet)$ has a left adjoint: the localization functor
$\Dcal_\lambda\otimes_{\Gamma(\Fl,\Dcal_\lambda)}\bullet$.
The Beilinson-Bernstein theorem states that the functor $\Gamma(\Fl,\bullet)$
is an equivalence provided $\lambda$ is regular and dominant, i.e., $\langle\lambda+\rho,\alpha\rangle\not\in \Z_{\leqslant 0}$ for any $\alpha\in \Delta$.

Let us explain some details on Ginzburg's localization theorem. For more details the reader is referred to \cite{Ginzburg_HC}.

One can consider the sheaf $\Dcal_\lambda$
as a quantization of the symplectic variety $T^*\Fl$. An analog of $T^*\Fl$ for $\Walg$ is the {\it Slodowy
variety} defined as follows. The action of $G$ on $T^*\Fl$ is Hamiltonian, the Springer resolution morphism
$\mu:T^*\Fl\rightarrow \g^*$ is a moment map. Recall the projection $\pi:\g^*\twoheadrightarrow \m^*$.
Then $\pi\circ\mu$ is a moment map for the $M$-action on $T^*\Fl$. By definition, the Slodowy variety $\mathcal{S}$
is the Hamiltonian reduction $(\pi\circ\mu)^{-1}(\chi|_{\m})/M$.  %There is a natural morphism
%$\mathcal{S}\rightarrow S$, its image is the intersection $S\cap\mathcal{N}$, where $\mathcal{N}$
%stands for the nilpotent cone in $\g^*$. One can show that $\mathcal{S}$ is a resolution of singularities
%of $S$.

To define an analog of the sheaf $\Dcal_\lambda$ in the W-algebra setting Ginzburg uses
the language of {\it directed algebras} (one can also use the language of microlocal sheaves,
see \cite{DK}). Once this analog is defined the Beilinson-Bernstein theorem transfers to the W-algebra
setting  verbatim. The scheme of the proof is as follows:  one introduces the notion of a Whittaker $\Dcal_\lambda$-module,
shows that the functors in the Beilinson-Bernstein theorem restrict to equivalences between the Whittaker
subcategories, and then uses the Skryabin theorem.

A related development is as follows.
Let $L$ be a finite dimensional $\g$-module, and $M$ be a Whittaker $\g$-module.
Then $L\otimes M$ is again a Whittaker $\g$-module. This allows to define tensor
products of finite dimensional $\g$-modules with $\Walg$-modules. These tensor products are studied
in detail in \cite{Goodwin}.

\subsection{Category $\Ocat$ for W-algebras}\label{SUBSECTION_Ocat}
In the representation  theory of $\U$ a crucial role is played by the category $\Ocat$
established by Bernstein, I. Gelfand and S. Gelfand in \cite{BGG}. There is an analog of
the BGG category $\Ocat$ for $\Walg$ introduced by Brundan, Goodwin and Kleshchev
in \cite{BGK}. The most important result about this category is that it is equivalent
to a certain category of {\it generalized Whittaker} $\U$-modules, \cite{LOCat}.
Our exposition follows \cite{LOCat}.

Recall the group $Q:=Z_G(e,h,f)$ acting on $\Walg$ and an embedding $\q\hookrightarrow \Walg$,
see Subsections \ref{SUBSECTION_Uge_remarks},\ref{SUBSECTION_equiv_Slodowy}.
Pick a Cartan subalgebra $\t\subset \q$ and set
$\lf:=\z_\g(\t)$. Then $\lf$ is a minimal Levi subalgebra in $\g$ containing $e$.
Further, pick an integral (=lying in the character lattice of the corresponding maximal
torus of $Q$) element $\theta\in \t$ with $\z_\g(\theta)=\lf$.
A category we are going to consider will depend on $\theta$.

Consider the decomposition $\Walg=\bigoplus_{i\in \Z}\Walg_i$, where $\Walg_i:=\{w\in \Walg| [\theta,w]=iw\}$. Set $\Walg_{\geqslant 0}:=\bigoplus_{i\geqslant 0}\Walg_i, \Walg_{>0}:=\bigoplus_{i>0}\Walg_i,
\Walg_{\geqslant 0}^+:=\Walg_{\geqslant 0}\cap \Walg\Walg_{>0}$. Then $\Walg_{\geqslant 0}$ is a subalgebra
in $\Walg$, while $\Walg_{>0}$ and $\Walg_{\geqslant 0}^+$ are two-sided ideals in $\Walg_{\geqslant 0}$.

We say that a $\Walg$-module $N$ belongs to the category $\Ocat(\theta)$ (in \cite{LOCat} this category was denoted
by $\widetilde{\Ocat}^{\t}(\theta)$) if
\begin{itemize}
\item $N$ is finitely generated.
\item $\t\subset \Walg$ acts on $N$ by diagonalizable endomorphisms.
\item $\Walg_{>0}$ acts on $N$ by locally nilpotent endomorphisms.
\end{itemize}

For example, when $e$ is distinguished (i.e., $\q=\{0\}$), then $\Ocat(\theta)$ consists precisely
of all finite dimensional $\Walg$-modules. In this case the
notion of the category $\Ocat$ is pretty useless. The other extreme is the case when
$e$ is principal in $\lf$. We will see below that in this case we can say a lot about $\Ocat(\theta)$.

Let us present an important construction of a module in $\Ocat(\theta)$.
Pick a $\Walg_{\geqslant 0}/\Walg_{\geqslant 0}^+$-module $N^0$ with diagonalizable $\t$-action (e.g., irreducible).  Define the {\it Verma} module $\Delta^\theta(N^0)$ by $\Delta^\theta(N^0):=\Walg\otimes_{\Walg_{\geqslant 0}} N^0$.

The properties of $\Ocat(\theta)$ are quite expectable.

\begin{Prop}\label{Prop:Ocat_prop}
\begin{enumerate}
\item If $N^0$ is irreducible, then $\Delta^\theta(N^0)$ has a unique irreducible quotient, say
$L^\theta(N^0)$.
\item Any irreducible module in $\Ocat(\theta)$ is isomorphic to  $L^\theta(N^0)$
for unique $N^0$.
\item Let $N\in \Ocat(\theta)$
be such that all $\t$-eigenspaces in $N$ are finite dimensional. Then
$N$ has finite length.
\item  $\Delta^\theta(N^0)$ with $\dim N^0<\infty$ satisfies the assumptions of (3).
\end{enumerate}
\end{Prop}
This is proved in \cite{BGK}, Theorem 4.5, Corollary 4.12 (in \cite{BGK}  a bit different
definition was used, in particular, the assumption in (3) was a part of the definition,
but this does not matter).

The most crucial property of $\Ocat(\theta)$ is that it is equivalent to a certain category of $\U$-modules.
To define this category we need some more notation.

Let $\g=\bigoplus_{i\in \Z}\g_i$ be the decomposition into the eigenspaces of $\ad\theta$. In particular,
$\lf=\g_0$. Form the subalgebra $\underline{\m}\subset \g_0$ by analogy with $\m\subset \g$ but using
the pair $(\g_0,e)$ instead of $(\g,e)$.  We define the W-algebra $\Walg^0:=U(\g_0,e)$.
This notation is different from \cite{LOCat} but agrees with \cite{Miura}.
Consider the subalgebra $\widetilde{\m}:=\underline{\m}\oplus \g_{>0}\subset\g$ (where $\g_{>0}:=\bigoplus_{i>0}\g_i$)
and set $\widetilde{\m}_\chi:=\{\xi-\langle\chi,\xi\rangle, \xi\in \widetilde{\m}\}$. The element
$\chi\in \g^*$ is $\t$-invariant and so vanishes on $\g_{>0}$. Hence $\widetilde{\m}_\chi=\underline{\m}_\chi\times \g_{>0}$.

We say that a $\U$-module $M$ is a {\it generalized Whittaker} module (for $e$ and $\theta$) if
\begin{itemize}
\item $M$ is finitely generated.
\item $\t$ acts on $M$ by diagonalizable endomorphisms.
\item $\widetilde{\m}_\chi$ acts on $M$ by locally nilpotent endomorphisms.
\end{itemize}

The category of generalized Whittaker modules will be denoted by $\operatorname{Wh}(\theta)$
(this notation is again different from the one used in \cite{LOCat}).

Again, one can define   a Verma module in $\operatorname{Wh}(\theta)$. Let $N^0$ be a $\Walg^0$-module
with diagonalizable $\t$-action. Let $\Sk_0:\Walg^0$-$\operatorname{Mod}\rightarrow U(\g_0)$-$\operatorname{Mod}$ be the Skryabin functor (for the pair $\g_0,e$). Define the Verma module $\Delta^{e,\theta}(N^0):=\U\otimes_{U(\g_{\geqslant 0})}\Sk_0(N^0)$, where $U(\g_{\geqslant 0})$ acts on $\Sk_0(N^0)$ via a natural epimorphism
$U(\g_{\geqslant 0})\twoheadrightarrow U(\g_0)$.

The following theorem is (a part of) the main result of \cite{LOCat}.

\begin{Thm}[\cite{LOCat}, Theorem 4.1]\label{Thm_Ocat}
There is an isomorphism $\Psi:\Walg^0\rightarrow \Walg_{\geqslant 0}/\Walg_{\geqslant 0}^+$ and an equivalence
$\mathcal{K}: \operatorname{Wh}(\theta)\rightarrow \Ocat(\theta)$ of abelian categories
such that the functors $\mathcal{K}(\Delta^{e,\theta}(\bullet))$ and $\Delta^\theta(\Psi_*(\bullet))$ from the category
of $\t$-diagonalizable $\Walg^0$-modules to $\Ocat(\theta)$ are isomorphic.
Here $\Psi_*$ denotes the push-forward functor with respect to the isomorphism $\Psi$.
\end{Thm}

Let us make a remark on an isomorphism $\Psi$. Such an isomorphism was first established
in \cite{BGK}. It is however not completely clear if one can use the isomorphism from
\cite{BGK} in Theorem \ref{Thm_Ocat}. A peculiar feature of both isomorphisms
is that they do not intertwine the embeddings $\t\hookrightarrow \Walg_{\geqslant 0}/\Walg_{\geqslant 0}^+,\Walg^0$ but rather induce a shift on $\t$. Since this shift will be of importance later we will give some details, see
Remark 5.5 in \cite{LOCat}. Namely, let $\iota^0,\iota$ denote the embeddings of $\t$ to $\Walg^0,\Walg_{\geqslant 0}/\Walg_{\geqslant 0}^+$, respectively. Then we have $\iota(\xi)=\Psi(\iota^0(\xi))-\langle\delta,\xi\rangle$
for an element $\delta\in \t^*$ constructed as follows. Pick a Cartan subalgebra $\h\subset\g$
containing $\t$ and $h$. Let $\Delta_{<0}$ denote the set of all roots $\alpha$ of $\g$ with
$\langle\alpha,\theta\rangle<0$. Set
\begin{equation}\label{eq:delta}
\delta:= \sum_{\alpha\in \Delta_{<0},\langle \alpha,h\rangle=1}\frac{1}{2}\alpha|_\t+ \sum_{\alpha\in \Delta_{<0},\langle\alpha,h\rangle\geqslant 2}\alpha|_{\t}.
\end{equation}

Till the end of the subsection we consider the category $\operatorname{Wh}(\theta)$ in  the special case when $e$ is principal in $\lf$.
Here  $\operatorname{Wh}(\theta)$ (with a slightly different
definition) was studied before by McDowell, \cite{McD},  by Milicic and Soergel, \cite{MS}, and by Backelin, \cite{Backelin}.

%In the special case when $e$ is principal in $\lf$ Theorem \ref{Thm_Ocat}  together with result
%of Backelin, \cite{Backelin}, allows to compute the multiplicities in $\Ocat(\theta)$.
To proceed we need some more notation.
Choose a system $\Pi$ of simple roots such that $\theta$ is dominant.
Then $\Pi_0:=\{\alpha\in \Pi: \langle\alpha, \theta\rangle=0\}$ is a  system of simple roots for
$\lf$. Let $\Delta_+,\Delta_{\lf+}$ denote the systems of positive roots for $\g$ and $\lf$.
For a root $\alpha$ let $e_\alpha$ denote a corresponding weight vector in $\g$.
Further, let $W_{\lf}$ denote the Weyl group of $\lf$. We have the dot action of
$W$ on $\h^*$ defined as in the previous subsection.

The W-algebra $\Walg^0$ is identified with the center $Z(\lf)$ of
$U(\lf)$. So all irreducible $\Walg^0$-module are 1-dimensional. The set of their
isomorphism classes is in one-to-one correspondence with orbits of the dot action of
$W_{\lf}$ on $\h^*$.

One may assume that $e=\sum_{\alpha\in \Pi_0}e_{-\alpha}$. Then $\widetilde{\m}$ is nothing else
but the maximal nilpotent subalgebra $\n$ of $\g$ corresponding to $\Pi$. Also we have $\langle\chi,e_\alpha\rangle\neq 0$ if and only if $\alpha\in \Pi_0$. So we recover the setting of \cite{Backelin},\cite{McD},\cite{MS}.

For $\lambda\in \h^*$ let us write $\Delta(\lambda),L(\lambda)$ for the Verma and irreducible
modules with highest weight $\lambda$ in the BGG category $\Ocat$  and $\Delta^{e,\theta}(\lambda), L^{e,\theta}(\lambda)$ for the Verma and irreducible modules in $\operatorname{Wh}(\theta)$
corresponding to $W_{\lf}\cdot \lambda$.

In \cite{MS} Milicic and Soergel proved that the (infinitesimal) block of $\operatorname{Wh}(\theta)$ corresponding to a regular integral central character is equivalent to the block of the BGG category $\Ocat$
with certain {\it singular} integral character that can be recovered from $\Pi_0$.
For a generalization of this equivalence  to the parabolic setting see \cite{Webster}.

For other blocks in $\operatorname{Wh}(\theta)$ (corresponding to singular/non-integral central
characters) the situation is more subtle. But still one can relate the multiplicities
in $\Ocat$ and in $\operatorname{Wh}(\theta)$. For $\lambda,\mu\in \h^*$ let
$[\Delta(\lambda):L(\mu)], [\Delta^{e,\theta}(\lambda): L^{e,\theta}(\mu)]$ denote the multiplicities
in the corresponding categories.

\begin{Thm}[\cite{Backelin}, Theorem 6.2]\label{Thm:4.2}
Let $\lambda,\mu\in \h^*$. If
\begin{enumerate}\item $\lambda\in W\cdot \mu$,
\item  and there is $w\in W_{\lf}$ such that
 $w\cdot \mu$ is antidominant
for $\lf$ and $\lambda-w\cdot \mu\in \Span_{\Z_{\geqslant 0}}(\Delta^+)$,
\end{enumerate}
then $[\Delta^{e,\theta}(\lambda):
L^{e,\theta}(\mu)]=[\Delta(\lambda):L(w\cdot\mu)]$. Otherwise,
$[\Delta^{e,\theta}(\lambda): L^{e,\theta}(\mu)]=0$.
\end{Thm}

An element $\lambda\in\h^*$ is said to be antidominant for $\lf$ if $\langle \lambda+\rho,\alpha^\vee\rangle\not\in \Z_{>0}$ for any  $\alpha\in \Delta_{\lf+}$.

%\subsection{Category of generalized Whittaker modules}\label{SUBSECTION_gen_Whit}
%\subsection{Category equivalence  for $\Ocat$}\label{SUBSECTION_equiv_OCat}
%\subsection{Localization theorem}

\section{Ideals in $U(\g)$ versus ideals in $\Walg$}\label{SECTION_ideals}
In this section we will  construct maps between the sets $\Id(\U)$ and $\Id(\Walg)$ of two-sided ideals in $\U$ and $\Walg$, respectively. This is done in the first two subsections.
In Subsection \ref{SUBSECTION_classification} we explain how these maps allow to relate
 (isomorphism classes of) finite dimensional irreducible $\Walg$-modules to
 primitive ideals  $\J\subset\U$  such that the associated variety
$\VA(\U/\J)$ is $\overline{\Orb}$. We conclude the section with some remarks in Subsection \ref{SUBSECTION_ideals_remarks}.

\subsection{Map $\bullet_\dagger: \Id(\U)\rightarrow \Id(\Walg)$}
Recall the algebras $\U_\hbar,\W_\hbar$ from Subsection \ref{SUBSECTION_Fedosov}, $\Walg_\hbar$
 from Subsection \ref{SUBSECTION_equiv_Slodowy}, and  the topological algebras $\U^\wedge_\hbar, \W^\wedge_\hbar, \Walg^\wedge_\hbar$ from Subsection \ref{SUBSECTION_Decomp}. By Theorem  \ref{Thm:decomposition},
 $\U^\wedge_\hbar\cong \W^\wedge_\hbar\widehat{\otimes}_{\K[[\hbar]]}\Walg^\wedge_\hbar$.
 Let us introduce suitable sets of ideals
of $\U_\hbar, \U^\wedge_\hbar, \Walg^\wedge_\hbar,\Walg_\hbar$. Namely, let $\Id_\hbar(\U_\hbar)$
denote the set of $\K^\times$-stable $\hbar$-saturated two-sided ideals in $\U_\hbar$ (an ideal $\J_\hbar\subset\U_\hbar$ is said to be $\hbar$-saturated if $\hbar a\in \J_\hbar$ implies $a\in \J_\hbar$,
in other words, if the quotient $\U_\hbar/\J_\hbar$ is a flat $\K[\hbar]$-module).  Define
the sets $\Id_\hbar(\U^\wedge_\hbar),\Id_\hbar(\Walg_\hbar), \Id_\hbar(\Walg^\wedge_\hbar)$ in a similar way. We define a map $\bullet_\dagger$ as the composition
\begin{equation}\label{eq:composition}\Id(\U)\xrightarrow{(a)} \Id_\hbar(\U_\hbar)\xrightarrow{(b)} \Id_\hbar(\U^\wedge_\hbar)\xrightarrow{(c)}
\Id_\hbar(\Walg^\wedge_\hbar)\xrightarrow{(d)} \Id_\hbar(\Walg_\hbar)\xrightarrow{(e)}\Id(\Walg).\end{equation}
Let us describe the intermediate maps.

{\it (a):} this map sends $\J\in \Id(\U)$ to $R_\hbar(\J):=\bigoplus (\J\cap \F_i\U)\hbar^i$. It
is a bijection, the inverse map sends $\J_\hbar\in \Id_\hbar(\U_\hbar)$ to its image under the natural
epimorphism $\U_\hbar\twoheadrightarrow \U_\hbar/(\hbar-1)=\U$.

{\it (b):}  this map sends $\J_\hbar\in \Id_\hbar(\U_\hbar)$  to its closure $\J_\hbar^\wedge\subset \U^\wedge_\hbar$. Equivalently, $\J_\hbar^\wedge=\U^\wedge_\hbar\J_\hbar$. This map
is neither injective (but it is easy to say when two ideals have the same image) nor surjective.

{\it (c):} this map sends $\J'_\hbar\in \Id_\hbar(\U_\hbar^\wedge)$ to $\I'_\hbar:=\J'_\hbar\cap \Walg^\wedge_\hbar$.
 It is a bijection: its inverse sends  $\I'_\hbar\in \Id_\hbar(\Walg^\wedge_\hbar)$ to $\W_\hbar^\wedge\widehat{\otimes}_{\K[[\hbar]]}\I'_\hbar$.

{\it (d):} this map sends $\I'_\hbar\in \Id_\hbar(\Walg^\wedge_\hbar)$ to $\I_\hbar:=\I'_\hbar\cap \Walg_\hbar$.
It is again a bijection, its inverse sends $\I_\hbar\in \Id_\hbar(\Walg_\hbar)$ to its closure.

{\it (e):} this map is analogous to the inverse of (a).

%The properties of the map $\J\rightarrow \J_\dagger$ are listed in the following proposition.

\begin{Prop}\label{Prop_lower_dag}
The map $\J\mapsto\J_\dagger$ has the following properties.
\begin{itemize}
\item[(1)] $\J_\dagger$ is $Q$-stable.
\item[(2)] $\gr \Walg/\J_\dagger$ is the pull-back of the $\K[\g^*]$-module
$\gr\U/\J$ to $S\subset\g^*$.
\item[(3)] $\J_\dagger=\Walg$ if and only if $\Orb\cap \VA(\U/\J)=\varnothing$.
\item[(4)] $\J_\dagger$ is a proper ideal of finite codimension in $\Walg$ if and only if $\overline{\Orb}$ is an irreducible
component of $\VA(\U/\J)$. In this case $\dim\Walg/\J_\dagger$  equals the multiplicity
of $\U/\J$ on $\Orb$.
\item[(5)] The natural map $(\J/\J\m_\chi)^{\ad \m}\rightarrow (\U/\U\m_\chi)^{\ad \m}$ is injective. Its image coincides
with $\J_\dagger$.
\end{itemize}
\end{Prop}
(1) follows directly from the construction. (2) is Proposition 3.4.2 in \cite{Wquant}. (3) and (4)
follow from (2).  (5) follows from Subsection 3.5 in \cite{HC}.

\subsection{Map $\bullet^\dagger: \Id(\Walg)\rightarrow \Id(\U)$}
By definition, $\bullet^\dagger$
is the composition
$$\Id(\Walg)\rightarrow \Id_\hbar(\Walg_\hbar)\rightarrow \Id_\hbar(\Walg^\wedge_\hbar)\rightarrow
\Id_\hbar(\U^\wedge_\hbar)\rightarrow \Id_\hbar(\U_\hbar)\rightarrow \Id(\U),$$
where all  maps except $\Id_\hbar(\U^\wedge_\hbar)\rightarrow \Id_\hbar(\U_\hbar)$ are the inverses of the corresponding
maps in (\ref{eq:composition}). The map $\Id_\hbar(\U^\wedge_\hbar)\rightarrow \Id_\hbar(\U_\hbar)$ sends
$\J'_\hbar\in \Id_\hbar(\U^\wedge_\hbar)$ to $\J'_\hbar\cap \U_\hbar$.

Let us list some properties of the map $\I\mapsto \I^\dagger: \Id(\Walg)\rightarrow \Id(\U)$.

\begin{Prop}\label{Prop_upper_dag}
\begin{itemize}
\item[(1)] Let $N$ be a $\Walg$-module. Then $\Ann_\Walg(N)^\dagger= \Ann_\U(\Sk(N))$, where $\Sk$
denotes the Skryabin functor, see Subsection \ref{SUBSECTION_Skryabin}.
\item[(2)] Let $N$ be a $\Walg$-module from the category $\Ocat(\theta)$, see Subsection \ref{SUBSECTION_Ocat}.
Then $\Ann_{\Walg}(N)^\dagger=\Ann_U(\Kfun(N))$, where $\Kfun$ is the functor from
Theorem \ref{Thm_Ocat}.
\item[(3)] Let $\I$ be an ideal of finite codimension in $\Walg$. If $\I$ is prime (resp., completely prime, primitive), then so is $\I^\dagger$.
\item[(4)] $V(\U/\I^\dagger)=\overline{\Orb}$ if and only if $\I$ is of finite codimension.
\item[(5)] Recall that the center $\Centr$ of $\U$ is identified with the center of $\Walg$.
Under this identification for any $\I\in \Id(\Walg)$ we have $\I\cap \Centr=\I^\dagger\cap \Centr$.
\item[(6)] The map $\I\mapsto \I^\dagger$ is $Q$-invariant.
\end{itemize}
\end{Prop}
Recall that an ideal $I$ in an associative unital algebra $A$ is called prime (resp., completely prime)
if $a$ or $b$ lies in $A$ whenever $aAb\subset I$ (resp., $ab\in I$). An ideal $I$ is said to be primitive if it
is the annihilator of an irreducible $A$-module.

(1) is assertion (ii) of \cite{Wquant}, Theorem 1.2.2. (2) is a part of \cite{LOCat}, Theorem 4.1.
(3) stems from \cite{Wquant}, Theorem 1.2.2. The "if" part of (4) follows from (1) and \cite{Premet2},
Theorem 3.1. The "only if" part follows from the inclusion $(\I^\dagger)_\dagger\subset \I$
that is a  direct consequence of our constructions.
(5) is assertion (iii) of \cite{Wquant}, Theorem 1.2.2. (6) follows directly from the construction.

\subsection{Classification of finite dimensional irreducible $\Walg$-modules}\label{SUBSECTION_classification}
This subsection is perhaps the most important part of the notes. Here we explain known results
about the classification of finite dimensional irreducible $\Walg$-modules. We have two results here, both are due to the author,
\cite{HC},\cite{Miura}. Both relate the set $\Irr_{fin}(\Walg)$ of (isomorphism classes of) finite dimensional irreducible $\Walg$-modules
to the set $\Prim_\Orb(\U)$ consisting of all primitive ideals $\J\subset \U$ with $\VA(\U/\J)=\overline{\Orb}$.

The first result was conjectured by Premet (private communication). To state it we notice that  the set
$\Irr_{fin}(\Walg)$ is canonically identified with the set $\Prim_{fin}(\Walg)$ of maximal (=primitive) ideals of finite codimension
in $\Walg$ (via taking the annihilator).
 Thanks to assertions (3),(4)
of Proposition \ref{Prop_upper_dag}, we see that
$N \mapsto \Ann_\Walg(N)^\dagger$ is a map $\Irr_{fin}(\Walg)\rightarrow \Prim_{\Orb}(\U)$. The group
$Q$  acts on $\Irr_{fin}(\Walg)$. The connected component $Q^\circ$ of $Q$ acts trivially
because the corresponding action of $\q$ on $\Walg$ is by inner derivations. So  the  $Q$-action on $\Irr_{fin}(\Walg)$ descends to that of the component group $C(e)=Q/Q^\circ$. By assertion (6)
of Proposition \ref{Prop_upper_dag}, the map $\Irr_{fin}(\Walg)\rightarrow \Prim_\Orb(\U)$ is $C(e)$-invariant.

\begin{Conj}[Premet]\label{Conj:6}
The map $N\mapsto \Ann(N)^\dagger: \Irr_{fin}(\Walg)\rightarrow \Prim_{\Orb}(\U)$
is surjective and any of its fibers is a single $C(e)$-orbit.
\end{Conj}

In \cite{Premet3} Premet proved that any $\J\in \Prim_{\Orb}(\U)$ with {\it rational} central character
lies in the image.
In  full generality the surjectivity part was first
proved in \cite{Wquant}, Theorem 1.2.2. Later alternative proofs were found
in \cite{Ginzburg_HC},\cite{Premet4}.
The description of fibers is more subtle. It  was obtained in \cite{HC}. It is a corollary of the following theorem.

%Consider the subset $\Id_{\Orb}(\U)\subset \Id(\U)$ consisting of all ideals $\J$ with $\VA(\U/\J)=\overline{\Orb}$
%and the subset $\Id_{fin}(\Walg)\subset \Id(\Walg)$ consisting of all ideals of finite codimension.

\begin{Thm}[\cite{HC}, Theorem 1.2.2]\label{Thm_rel_btw_id}
Let $\I$ be a {\rm $Q$-stable} ideal of finite codimension in $\Walg$. Then
$(\I^\dagger)_\dagger=\I$.
\end{Thm}

The second  result is stated in terms of the category $\Ocat(\theta)$. Let $\theta,\lf=\g_0,\Walg^0$ have the same meaning as  in
Subsection \ref{SUBSECTION_Ocat}. Choose a Cartan subalgebra $\h\subset \lf$ and a system of simple roots
$\Pi\subset \h^*$ as in the discussion preceding Theorem \ref{Thm:4.2}.

Let us introduce some more notation. For $\lambda\in \h^*$ let  $L_0(\lambda)$
stand for  the irreducible  $\g_0$-module with highest weight $\lambda$.
Set $J(\lambda):=\Ann_\U(L(\lambda)), J_0(\lambda):=\Ann_{U(\g_0)}(L_0(\lambda))$. According to Duflo, \cite{Duflo}, any primitive ideal in $\U$ (resp., in $U(\g_0)$) has the form $J(\lambda)$ (resp., $J_0(\lambda)$)
for some (in general, non-unique)  $\lambda\in \h^*$.

\begin{Prop}\label{Prop:hw_dag}[\cite{Miura}, Theorem 5.1.1]
Let $N_0$ be an irreducible finite dimensional $\Walg^0$-module. If $\Ann_{\Walg^0}(N_0)^\dagger= J_0(\lambda)$ for some $\lambda$, then $\Ann_\Walg(L^\theta(N_0))^\dagger= J(\lambda)$. In particular,
$L^\theta(N_0)$ is finite dimensional if and only if $\VA(\U/J(\lambda))=\overline{\Orb}$.
\end{Prop}

\subsection{Remarks}\label{SUBSECTION_ideals_remarks}
In the representation theory of $\U$ there are many results on the computation
of $\VA(\U/J(\lambda))$ and on the description of $\Prim_\Orb(\U)$. They are due to Joseph, Barbasch-Vogan and
others, see, for example, \cite{BV1},\cite{BV2},\cite{Joseph}.
In particular, it is known that $\Prim_\Orb(\U)$ is always non-empty.

Next, we remark that the maps between the sets of ideals upgrade to functors between the categories
of Harish-Chandra bimodules, see \cite{Ginzburg_HC},\cite{HC}. The study of these functors is supposed to help
to obtain the complete description of $\Irr_{fin}(\Walg)$ itself (not just of the quotient
$\Irr_{fin}(\Walg)/C(e)$).

\section{One-dimensional $\Walg$-modules}\label{SECTION_one_dim}
\subsection{Motivation}
The following conjecture was made by Premet.
\begin{Conj}[\cite{Premet2}, Conjecture 3.1]\label{Conj:7.0}
Any W-algebra has a one-dimensional representation (equivalently, a two-sided ideal
of codimension 1).
\end{Conj}

At the moment when this text is being written the conjecture is known to be true with exception of
several cases in type $E_8$, where it is still open.

The reason why Conjecture \ref{Conj:7.0} is important is that it implies affirmative answers to some
old questions in representation theory of universal enveloping algebras:
\begin{itemize}
\item[(A)] the  question of Humphreys on the   existence of a small non-restricted representation
for semisimple Lie algebras in  characteristic $p$.
\item[(B)]
the existence of a completely
prime primitive ideal with given associated variety (this question traces back, at least, to Dixmier)
\end{itemize}

The proof that Conjecture \ref{Conj:7.0} implies
(A) for $p\gg 0$  is obtained in \cite{Premet4}, Theorem 1.4.

The claim that Conjecture \ref{Conj:7.0} implies (B) follows from Proposition \ref{Prop_upper_dag}:
if $\I\subset \Walg$ has codimension $1$, then $\I^\dagger$ is completely prime and $\VA(\U/\I^\dagger)=\overline{\Orb}$.

In fact, the implication in the previous paragraph was obtained earlier by Moeglin, \cite{Moeglin1},\cite{Moeglin2}.
She considered primitive ideals in $\U$ admitting a {\it Whittaker model}. Using the techniques of \cite{GG},
one can show that a Whittaker model in the sense of Moeglin is precisely the image of a one-dimensional
$\Walg$-module under the Skryabin equivalence.

Actually, Moeglin obtained a  stronger result: that any ideal admitting a Whittaker model gives rise to
a unique quantization (in an appropriate sense) of a suitable covering of $\Orb$, see \cite{Moeglin2}
for details.

\subsection{Classical algebras}
It turns out that Conjecture \ref{Conj:7.0} holds for all classical simple Lie algebras.
This was proved in \cite{Wquant}, Theorem 1.2.3, (1). Let us describe the idea of the proof.

We need to show that there is an ideal of codimension 1 in $\Walg$. Thanks to Proposition \ref{Prop_lower_dag},
this is the case when there is $\J\in \Id(\U)$  such that $\overline{\Orb}$ is an irreducible component
of $\VA(\U/\J)$ and the multiplicity of $\U/\J$ on $\overline{\Orb}$ is 1 (this implication
also was proved by Moeglin using the language of Whittaker models, see \cite{Moeglin2}).

Let $G$ be one of the classical groups $\operatorname{SL}_n(\K),$ $\operatorname{O}_n(\K), \operatorname{Sp}_{2n}(\K)$ (depending on $\g$). We emphasize that for $\g=\mathfrak{so}_n$ we need a disconnected group. It turns out that there is an ideal $\J$ in $\U$ such that $\gr \U/\J=\K[\overline{Ge}]$,
where $\gr$ is taken with respect to the filtration on $\U/\J$ induced from the PBW filtration on $\U$.
Such an ideal is obtained by the quantization of the Kraft-Procesi construction of $\overline{Ge}$
via a Hamiltonian reduction of a vector space, see \cite{KP1},\cite{KP2}\footnote{After \cite{Wquant} was already published I learned that the construction of $\J$ used there (and explained above) was discovered before by R. Brylinski,\cite{Brylinski}.}.

In type $A$ more can be said. Form the quotient
$\Walg^{ab}$ of $\Walg$ by the relations $[a,b], a,b\in \Walg$. The one-dimensional $\Walg$-modules are parametrized by
points of $\Spec(\Walg^{ab})$. In \cite{Premet4}, Subsection 3.8, Premet proved that
for $\g=\sl_n$ the algebra $\Walg^{ab}$ is the polynomial algebra in $d-1$ variables, where $d$ is the maximal
size of a Jordan block of $e$.
Premet's proof is based on the Brundan-Kleshchev presentation of $\Walg$, see Subsection \ref{SUBSECTION_A_Yangians} for details.

\subsection{Parabolic induction}
It is easy to prove Conjecture \ref{Conj:7.0} when $e$ is even. Indeed, as we have seen in Subsection \ref{SUBSECTION_Uge_even},
the algebra $\Walg$ for even $e$ can be embedded
into $U(\g(0))$, see Subsection  \ref{SUBSECTION_Uge_even}. Then we can take any 1-dimensional representation
of $U(\g(0))$ and restrict it to $\Walg$.

Premet, \cite{Premet4}, observed that a similar result holds in a much more general setting.
In the theory of nilpotent elements in semisimple Lie algebras there is a construction called the {\it Lusztig-Spaltenstein induction}. It was  introduced in \cite{LS}, for a review see, for example, \cite{McG}.
Namely, let $\underline{\g}\subset \g$ be a Levi subalgebra and $\underline{\Orb}\subset \underline{\g}$
be a nilpotent orbit. The Lusztig-Spaltenstein induction produces a nilpotent orbit
$\Orb\subset \g$ from the pair $(\underline{\g},\underline{\Orb})$. We say that $\Orb$ is induced
from $(\underline{\g},\underline{\Orb})$. If $e$ is even, then $\Orb$ is induced
 from $(\g(0),\{0\})$.  If $\Orb$ cannot be induced from a nilpotent orbit in proper Levi subalgebra,
$\Orb$ is called {\it rigid}.

\begin{Thm}[Premet, \cite{Premet4}]\label{Thm:7.1}
Let $\Orb$ be induced from $(\underline{\g},\underline{\Orb})$. If the algebra
$\underline{\Walg}:=U(\underline{\g},\underline{e})$, where $\underline{e}\in \underline{\Orb}$, has a one-dimensional representation, then  $\Walg$ does.
\end{Thm}

Premet's proof of Theorem \ref{Thm:7.1} is based on the
reduction to positive characteristic. Another proof, close in spirit to that for even elements, was found
by the author in \cite{Miura}. Namely, under the assumptions of Theorem \ref{Thm:7.1}, there is an embedding
of $\Walg$ into a certain {\it completion} of $\underline{\Walg}$. The latter acts on all finite dimensional
$\underline{\Walg}$-modules. So a one-dimensional $\Walg$-module again can be obtained by restriction.

\subsection{1-dimensional representations via category $\Ocat(\theta)$}
In this subsection we will explain how to apply the category $\Ocat(\theta)$ to the study
of one-dimensional representations of $\Walg$, see \cite{Miura}.

We use the notation from Subsection \ref{SUBSECTION_Ocat}.
Let us impose the following condition on a nilpotent element
$e$:
\begin{itemize}
\item[(*)] the algebra $\q$ is semisimple.
\end{itemize}
It turns out that this condition is satisfied for all rigid nilpotent elements. A proof based
on the classification of such elements can be found in \cite{Miura}, Subsection 5.2. It would be very
interesting to find a conceptual proof.

Let $N^0$ be a  finite dimensional
$\Walg^0$-module. We want a criterium for $L^\theta(N^0)$ to be 1-dimensional. Since
 $N^0\hookrightarrow L^\theta(N^0)$, of course,
$\dim L^\theta(N^0)=1$ implies $\dim N^0=1$.

The following result  follows from Theorem 5.2.1 in \cite{Miura}.

\begin{Thm}\label{Thm:one_dim_Ocat}
Suppose the condition (*) holds. Let $N^0$ be a 1-dimensional $\Walg^0$-module.
Then the following conditions are equivalent:
\begin{enumerate}
\item $\dim L^\theta(N^0)=1$.
\item $\t\subset \Walg^0$ acts on $N^0$ by $\delta$, where  $\delta$ is given by
(\ref{eq:delta}).
\end{enumerate}
\end{Thm}

When $e$ is of principal Levi type (which is true for all but 2 rigid nilpotent
elements in exceptional Lie algebras), then any irreducible $\Walg^0$-module
is 1-dimensional (recall that $\Walg^0$ is just the center of $U(\lf)$).

Combining Theorem \ref{Thm:one_dim_Ocat}, Proposition \ref{Prop:hw_dag}, and assertion (4) of Proposition \ref{Prop_upper_dag} one obtains a criterium for an ideal $\J\subset \U$ to have the
form $\I^\dagger$ with $\dim \Walg/\I=1$. More precisely, we have the following result,
\cite{Miura}, Subsection 5.3.

\begin{Cor}\label{Cor:1dim}
Suppose $\q$ satisfies (*). Let
$\theta,\h,\Pi$ be chosen as in the discussion preceding Theorem \ref{Thm:4.2}.
Let $\Orb_0$ denote the adjoint orbit of $e$ in $\lf$.
\begin{enumerate}
\item Let $\lambda\in \h^*$ satisfy the following four conditions:
\begin{itemize}
\item[(A)] The associated variety of $U(\lf)/J_0(\lambda)$ in $\g_0^*$ is $\overline{\Orb_0}$.
\item[(B)] $\dim \VA(\U/J(\lambda))\leqslant \dim \Orb$.
\item[(C)] $\lambda-\delta$ vanishes on $\t$.
\item[(D)] $J_0(\lambda)$ corresponds to an ideal of codimension 1 in $\Walg^0$.
\end{itemize}
Then $J(\lambda)=\I^\dagger$ for some ideal $\I\subset\Walg$ of codimension 1.
\item For any ideal $\I\subset \Walg$ of codimension 1 there is $\lambda\in \h^*$ satisfying
(A)-(D) and such that $J(\lambda)=\I^\dagger$.
\end{enumerate}
\end{Cor}
When $e$ is principal in $\lf$ the condition (A) means that $\lambda$ is antidominant for $\lf$,
while the condition (D) becomes vacuous.

\subsection{Exceptional algebras}
Let us summarize what is known about Conjecture \ref{Conj:7.0} for exceptional Lie algebras.
As Premet checked in \cite{Premet2}, $\Walg$ has a one-dimensional module
provided $e$ is a minimal nilpotent element (in an arbitrary simple Lie algebra).
His approach was based on an analysis of generators and relations for $\Walg$
that are not very difficult for minimal nilpotents. Recently Goodwin, R\"{o}hrle and
Ubly, \cite{GRU}, extended Premet's approach to all rigid nilpotents in $G_2,F_4,E_6,E_7$
and some rigid nilpotents in $E_8$. The result is that
for all nilpotent elements they considered a one-dimensional $\Walg$-module does exist. They used
the GAP program to analyze the relations. "Large" nilpotent elements in $E_8$ remain to complicated computationally.
Maybe, one can deduce  Conjecture \ref{Conj:7.0} for $E_8$ from Corollary \ref{Cor:1dim}.

\section{Type $A$}\label{SECTION_type_A}
This section is devoted to results concerning W-algebras for $\g=\sl_N$ (or $\g=\gl_N$).
In Subsection \ref{SUBSECTION_A_Yangians} we very briefly sketch a relation between
W-algebras and Yangians. In Subsection \ref{SUBSECTION_A_other} we mention some other results:
the higher level Schur-Weyl duality of Brundan and Kleshchev and the Gelfand-Kirillov conjecture
for W-algebras proved by Futorny, Molev and Ovsienko.
\subsection{W-algebras vs Yangians}\label{SUBSECTION_A_Yangians}
In this subsection we will briefly explain a relationship between $W$-algebras for $\g=\gl_N$
and certain inifinite dimensional algebras called {\it shifted Yangians}.
A shifted Yangian is a certain generalization of the usual Yangian for $\gl_n$.
For a comprehensive treatment of Yangians and related algebras the reader is referred to Molev's book
\cite{Molev}.   A relation between Yangians and W-algebras was first observed by Ragoucy and Sorba in \cite{RS}
and then generalized to shifted Yangians by Brundan and Kleshchev, \cite{BK1}.

The Yangian $Y(\gl_n)$ can be defined as the algebra generated by elements $t_{ij}^{(r)}, i,j=1,\ldots,n,
r\in \N,$ subject to the relations
$$[t_{ij}^{(r+1)},t^{(s)}_{kl}]-[t_{ij}^{(r)},t_{kl}^{(s+1)}]=t_{kj}^{(r)}t_{il}^{(s)}-t_{kj}^{(s)}t_{il}^{(r)}.$$

However, the generators $t_{ij}^{(r)}$ are not convenient to establish a relation between the Yangians and W-algebras.
In  \cite{BK_Yangian} Brundan and Kleshchev found a new presentation of $Y(\gl_n)$. Generalizing this presentation
they introduced shifted Yangians in \cite{BK1}.

A shifted Yangian $Y_n(\sigma)$ depends on a positive integer $n$ and some {\it shift matrix} $\sigma$.
By definition, $\sigma=(s_{ij})_{i,j=1}^n$ is a shift matrix if $s_{ij}$ is a nonnegative integer
("shift") with $s_{ij}+s_{jk}=s_{ik}$ whenever $|i-j|+|j-k|=|i-j|$. By definition, the algebra
$Y_n(\sigma)$ is given by generators $$D_{i}^{(r)} (1\leqslant i\leqslant n, r>0), E_{i}^{(r)} (1\leqslant i< n, r>s_{i,i+1}), F_{i}^{(r)} (1\leqslant i< n, r>s_{i+1,i})$$ subject to certain explicit relations
(see \cite{BK1}, (2.4)-(2.15)).  The shifted Yangian coincides with the usual one
when $\sigma=0$.
For $l>s_{1,n}+s_{n,1}$ define the quotient (the truncated shifted Yangian of level $l$) $Y_{n,l}(\sigma)$ of $Y_n(\sigma)$ by the two-sided ideal generated by $D_{i}^{(r)}, 1\leqslant i\leqslant n, r>p_1:=l-s_{1,n}-s_{n,1}$.

%Pick a Young diagram $\lambda = (\lambda_1\geq\cdots\geq\lambda_n)
%(one can also work with more general diagrams called pyramids...),
%and set l := \lambda_1. Then to n and \lambda one assigns the shift
%matrix \sigma^\lambda = (s_{ij})_{i,j=1}^n by setting s_{ij} := 0 for
%i >= j and s_{ij} := \lambda_{n+1-i} - \lambda_{n+1-j} for i <= j. In
%particular, for the Young diagram of shape n x l, we get \sigma^
%\lambda = 0.

To establish a relationship between shifted Yangians and  $W$-algebras fix a positive integer $n$, pick a Young diagram $\lambda=(\lambda_1,\ldots,\lambda_n), \lambda_1\geqslant \ldots\geqslant \lambda_n\geqslant 0$
(one can also work with more general diagrams called {\it pyramids},
see \cite{BK1}, $\S$7 for details), and set $l:=\lambda_1$.
Then to $n$ and $\lambda$ one can assign the {\it shift matrix} $\sigma^\lambda=(s_{ij})_{i,j=1}^n$
by setting $s_{ij} := 0$ for
$i \geqslant j$ and $s_{ij} := \lambda_{n+1-j}-\lambda_{n+1-i}$ for $i<j$.
In particular, for the Young diagram of shape $n\times l$, we get $\sigma=0$.

To $\lambda$ one assigns a nilpotent element $e_\lambda\in \gl_{N}$, where $N:=\sum_{i=1}^n \lambda_i$, in the usual way ($\lambda_i$ are the sizes of the Jordan blocks of $e_\lambda$).

\begin{Thm}[\cite{BK1}, Theorem 10.1]\label{Thm:Yangian}
$U(\gl_N,e_\lambda)\cong Y_{n,l}(\sigma^\lambda)$.
\end{Thm}

In \cite{BK2} Brundan and Kleshchev used this theorem to  study  the representation theory
of $U(\gl_N,e_\lambda)$. In particular, they obtained a classification of finite dimensional
irreducible $U(\gl_N,e)$-modules
(which also follows from Proposition \ref{Prop:hw_dag} thanks to Joseph's
computation of $\VA(U(\gl_N)/J(\lambda))$, see \cite{Joseph}; we remark that any nilpotent element
in $\gl_N$ is of principal Levi type).
%They also obtained the multiplicity
%formulas for the categories $\Ocat(\theta)$ for $U(\gl_n,e_\lambda)$, see Theorem \ref{Thm:4.2}.

There is a generalization of the results explained above in this subsection
to other classical Lie algebras first observed by Ragoucy, \cite{R}   and
worked out in more detail by
J. Brown, \cite{Brown1},\cite{Brown2}. Namely, for orthogonal and symplectic algebras
there are analogs of $Y(\gl_n)$ called {\it twisted} Yangians. Theorem \ref{Thm:Yangian} generalizes
to twisted Yangians. It is interesting that, similarly to
$Y(\gl_n)$-case, nilpotent elements arising in this generalization again correspond
to partitions with all parts equal. It is unclear whether there is a reasonable shifted
version of the twisted Yangians that is related to the W-algebras constructed from  arbitrary
nilpotent elements.

\subsection{Other results}\label{SUBSECTION_A_other}
W-algebras in type A enjoy some other interesting properties.

For example, in \cite{BK3} Brundan and Kleshchev obtained a very nice result: a "higher level" generalization
of the classical Schur-Weyl duality.
Recall that the classical Schur-Weyl duality  relates
between polynomial representations of $\operatorname{GL}_N(\K)$ and representations
of the symmetric group $S_d$ in $d$ letters. The Brundan-Kleshchev generalization
relates  modules over the cyclotomic degenerate Hecke algebra $H_{d}(\lambda)$ corresponding
to a partition $\lambda$ of $N$ (this algebra is a higher level generalization of $S_d$) and modules
over the W-algebra $U(\gl_N,e_\lambda)$. For details the reader is referred to \cite{BK3}
or to the review \cite{Wang} by Wang.

%The Brundan-Kleshchev duality can be applied to study the representation theory
%of

%It is unclear whether there is a reasonable generalization of the higher level Schur-Weyl duality
%for other types.

Another result we would like to mention is an analog of the Gelfand-Kirillov conjecture
for W-algebras proved in \cite{FMO}.

For a Noetherian domain $A$ let $Q(A)$ denote its skew-field of fractions.
Gelfand and Kirillov, \cite{GK}, conjectured   that for any finite dimensional
 algebraic Lie algebra $\mathfrak{a}$ the skew-field $Q(U(\mathfrak{a}))$ is isomorphic to $Q(\W_{l}(F_d))$,
where  $F_d$ is a purely transcendental extension of $\K$
of some degree $d$ and $\W_l(F_d)$ stands for the Weyl algebra of a $2l$-dimensional symplectic
vector space over $F_d$. In \cite{GK}  the conjecture was verified for $\g=\sl_n$\footnote{Recently Premet proved in \cite{Premet5}
that the Gelfand-Kirillov conjecture \underline{does not} hold for $\mathfrak{a}$ of type $B_n (n\geqslant 3),D_n (n\geqslant 4),E_6,E_7,E_8,F_4$.}.
In \cite{FMO} Futorny, Molev and Ovsienko proved that the straightforward analog of the Gelfand-Kirillov
conjecture holds for $U(\gl_n,e)$ (and for $U(\sl_n,e)$) for an arbitrary nilpotent element $e\in \sl_n$.

% This is for the references %
% The reference below to Ahlfor's Complex Analysis book %
% is a sample that you should remove and put your own   %
% references %
% You can cite it with \cite{ahlfors} %

% Done :-) %

\begin{thebibliography}{a}
\bibitem{Backelin} E. Backelin,  \emph{Representation of the category $\mathcal{O}$ in Whittaker
categories}, IMRN, 4(1997), 153-172.
\bibitem{BV1} D. Barbasch, D. Vogan, \emph{Primitive ideals and orbital integrals
in complex classical groups}, Math. Ann. 259(1982), 153-199.
\bibitem{BV2} D. Barbasch, D. Vogan, \emph{Primitive ideals and orbital integrals
in complex exceptional groups}, J. Algebra 80(1983), 350-382.
\bibitem{BB} A. Beilinson, J. Bernstein, \emph{Localization de
$\mathfrak{g}$-modules}, C. R. Acad. Sci. Paris 292 (1981), no. 1, 15-18.
\bibitem{BGG} I. Bernstein, I. Gelfand, S. Gelfand, \emph{A category of $\mathfrak{g}$-modules},
Funct. Anal. Appl. 10(1976), 87-92.
\bibitem{BT} J. de Boer, T. Tjin, \emph{Quantization and representation theory
of finite W-algebras}, Comm. Math. Phys. 158(1993), 485-516.
\bibitem{Brown1} J. Brown, \emph{Twisted Yangians and finite W-algebras}, Transform. Groups 14
(2009), 87-114.
\bibitem{Brown2} J. Brown, \emph{Representation theory of rectangular finite W-algebras},
arXiv:1003.2179.
\bibitem{BG} J. Brundan, S. Goodwin, \emph{Good gradings polytopes},
Proc. London Math. Soc. 94(2007), 155-180.
\bibitem{BGK} J. Brundan, S. Goodwin, A. Kleshchev, \emph{Highest weight theory
for finite $W$-algebras},  IMRN 2008, no. 15, Art. ID rnn051.
\bibitem{BK_Yangian} J. Brundan, A. Kleshchev, \emph{Parabolic presentation of the Yangian
$Y(\mathfrak{gl}_n)$}, Comm. Math. Phys. 254(2005), 191-220.
\bibitem{BK1} J. Brundan, A. Kleshchev, \emph{Shifted Yangians and finite $W$-algebras}, Adv. Math. 200(2006), 136-195.
\bibitem{BK2} J. Brundan, A. Kleshchev, \emph{Representations of shifted Yangians and finite W-algebras},
Mem. Amer. Math. Soc. 196 (2008), 107 pp.
\bibitem{BK3} J. Brundan, A. Kleshchev, \emph{Schur-Weyl duality
for higher levels}, Selecta Math., 14(2008), 1-57.
\bibitem{Brylinski} R. Brylinski, \emph{Dixmier algebras for classical
complex nilpotent orbits via Kraft-Procesi models. I},  Prog. Math. 213,
Birkh\"{a}user, Boston, 49-67.
\bibitem{DSK_appendix} A. D'Andrea, C. De Concini, A. De
Sole, R. Heluani and V. Kac, \emph{Three equivalent definitions of finite W-algebras},
Appendix to \cite{DSK}.
\bibitem{DSK} A. De Sole, V. Kac, \emph{Finite vs affine
$W$-algebras},  Japan. J. Math, 1(2006), 137-261.
\bibitem{DK} C. Dodd, K. Kremnizer, \emph{A Localization Theorem for Finite W-algebras},
arXiv:0911.2210.
\bibitem{Duflo} M. Duflo, \emph{Sur la classification des id\'{e}aux primitifs
dans l'alg\`{e}bre envellopante d'une alg\`{e}bre de Lie semi-simple}, Ann. Math. 105(1977), 107-120.
\bibitem{EK} A. Elashvili, V. Kac, \emph{Classification of good gradings of simple
Lie algebras}, in: "Lie groups and invariant theory" (E.B. Vinberg ed.), Amer. Math. Soc. Transl.
ser. 2, 213(2005), 85-104.
\bibitem{Fedosov1} B. Fedosov, \emph{A simple geometrical construction of
deformation quantization}, J. Diff. Geom. 40(1994), 213-238.
\bibitem{Fedosov2} B. Fedosov, \emph{Deformation quantization and index theory},
in Mathematical Topics 9, Akademie Verlag, 1996.
\bibitem{FMO} V. Futorny, A. Molev, S. Ovsienko, \emph{Gelfand-Kirillov conjecture
and Gelfand-Tsetlin modules for finite W-algebras}, Adv. Math. 223(2010), 773-796.
\bibitem{GK} I. Gelfand, A. Kirillov, \emph{Sur les corps l\'{e}s  aux alg\`{e}bres enveloppantes des alg\`{e}bres
de Lie}, Publ. IHES, 31(1966), 5-19.
\bibitem{GG} W.L. Gan, V. Ginzburg, \emph{Quantization of Slodowy slices}, IMRN, 5(2002), 243-255.
\bibitem{Ginzburg_HC} V. Ginzburg, \emph{Harish-Chandra bimodules for quantized Slodowy
slices}, Repres. Theory 13(2009), 236-271.
\bibitem{Goodwin} S. Goodwin, \emph{Translation for finite W-algebras}, arXiv:0908.2739.
\bibitem{GRU} S. Goodwin, G. R\"{o}hrle, G. Ubly, \emph{On 1-dimensional representations of finite W-algebras associated to simple Lie algebras of exceptional type}, arXiv:0905.3714.
\bibitem{Joseph} A. Joseph, \emph{Sur la classification des id\'{e}aux primitifs dans l'algebre envellopante
de $\sl(n+1,\mathbb{C})$}, C.R. Acad. Sci. Paris S\'{e}r A-B, 287(1978), N5, A303-306.
\bibitem{Kawanaka1} N. Kawanaka, \emph{Generalized Gelfand-Graev representations and Ennola duality}, In
"Algebraic Groups and Related Topics", Advanced Studies in Pure Mathematics 6(1985), North-Holland, p. 175-206.
\bibitem{Kostant_pol} B. Kostant, \emph{Lie group representations on polynomial rings},
Amer. J. Math. 85(1963), 327-404.
\bibitem{Kostant} B. Kostant, \emph{On Whittaker vectors and representation
theory}, Invent. Math. 48(1978), 101-184.
\bibitem{Kostant_Toda} B. Kostant,
\emph{The Solution to a Generalized Toda Lattice and Representation Theory}, Adv. in Math., 34(1979), 195-338.
\bibitem{KP1} H. Kraft, C. Procesi, \emph{Closures of conjugacy classes of matrices
are normal}, Invent. Math. 53 (1979), 227-247.
\bibitem{KP2} H. Kraft, C. Procesi, \emph{On the geometry of conjugacy classes
in classical groups}, Comment. Math. Helv. 57 (1982), 539-602.
\bibitem{slice} I.V. Losev, \emph{Symplectic slices for reductive groups},
Mat. Sbornik 197(2006), N2, 75-86 (in Russian). English translation in:
Sbornik Math. 197(2006), N2, 213-224.
\bibitem{Wquant} I.V. Losev, \emph{Quantized symplectic actions and
$W$-algebras}, J. Amer. Math. Soc. 23(2010), 35-59.
\bibitem{HC} I. Losev, \emph{Finite dimensional representations of W-algebras},
arXiv:0807.1023.
\bibitem{LOCat} I. Losev. \emph{On the structure of the category $\mathcal{O}$ for $W$-algebras},
arXiv:0812.1584.
\bibitem{Miura} I. Losev, \emph{1-dimensional representations and parabolic induction
for $W$-algebras}, arXiv:0906.0157.
%\bibitem{ES_appendix} I. Losev, Appendix to: P. Etingof, T. Schedler. {\it Poisson traces and D-modules on Poisson %varieties}, arXiv:0908.3868.
\bibitem{LS} G. Lusztig, N. Spaltenstein, \emph{Induced unipotent
classes}, J. London Math. Soc. (2), 19(1979), 41-52.
\bibitem{Lynch} T.E. Lynch, \emph{Generalized Whittaker vectors and representation
theory}, Thesis, M.I.T., 1979.
%\bibitem[Ma]{Matumoto} H. Matumoto. {\it Whittaker modules associated with hi}
\bibitem{McD} E. McDowell, \emph{On modules induced from Whittaker modules}, J. Algebra 96(1985), n.1, 161-177.
\bibitem{McG} W. McGovern, \emph{The adjoint representation and the
adjoint action}, Encyclopaedia of mathematical sciences, 131.
Invariant theory and algebraic transformation groups, II, Springer
Verlag, Berlin, 2002.
\bibitem{MS} D. Milicic, W. Soergel. \emph{The composition series of modules induced from
Whittaker modules}, Comment. Math. Helv. 72(1997), 503-520.
\bibitem{Moeglin1} C. Moeglin, \emph{Mod\`{e}les de Whittaker et id\'{e}aux primitifs compl\`{e}tement premiers dans les alg\`{e}bres enveloppantes I}, C.R. Acad. Sci. Paris, S\'{e}r. I 303(1986), No. 17,  845-848.
\bibitem{Moeglin2} C. Moeglin, \emph{Mod\`{e}les de Whittaker et id\'{e}aux primitifs compl\`{e}tement premiers dans les alg\`{e}bres enveloppantes II},  Math. Scand. 63 (1988), 5-35.
\bibitem{Molev} A. Molev, \emph{Yangians and classical Lie algebras}, Mathematical Surveys and Monographs, 142. AMS (2007).
\bibitem{Premet0} A. Premet, \emph{Irreducible representations of Lie algebras of reductive groups and the Kac-Weisfeiler conjecture}, Invent. Math. 121(1995), 79-117.
\bibitem{Premet1} A. Premet, \emph{Special transverse slices and their enveloping algebras}, Adv. Math. 170(2002),
1-55.
\bibitem{Premet2} A. Premet, \emph{Enveloping algebras of Slodowy slices and the Joseph ideal},
J. Eur. Math. Soc, 9(2007), N3, 487-543.
\bibitem{Premet3} A. Premet, \emph{Primitive ideals, non-restricted representations
and finite $W$-algebras}. Moscow Math. J.  7(2007), 743-762.
\bibitem{Premet4} A. Premet, \emph{Commutative quotients of finite W-algebras}, arXiv:0809.0663.
Accepted by Adv. Math.
\bibitem{Premet5} A. Premet, \emph{Modular Lie algebras and Gelfand-Kirillov conjecture}, arXiv:0907.2500.
\bibitem{R} E. Ragoucy, \emph{Twisted Yangians and folded W-algebras},
Internat. J. Modern. Phys. A 16 (2001), 2411-2433
\bibitem{RS} E. Ragoucy, P. Sorba, \emph{Yangian realizations from finite W-algebras},
Comm. Math. Phys. 203(1999), 551-572.
\bibitem{Skryabin} S. Skryabin, An appendix to \cite{Premet1}.
\bibitem{Slodowy} P. Slodowy, \emph{Simple singularities and simple algebraic
groups}, Lect. Notes Math., v.815. Springer, Berlin/Heidelberg/New
York, 1980.
\bibitem{VD} K. de Vos, P. van Driel, \emph{Kazhdan-Lusztig conjecture for finite W-algebras},
Lett. Math. Phys. 35(1996), 333-344.
\bibitem{Wang} W. Wang, \emph{Nilpotent orbits and W-algebras}, arXiv:0912.0689.
\bibitem{Webster} B. Webster, \emph{Singular blocks of parabolic category O and finite W-algebras},
arXiv:0909.1860.
\end{thebibliography}
\end{document}